\documentclass[twoside,reqno,a4paper]{amsart}
\usepackage{amsmath}
\usepackage{style}
\renewcommand{\P}{\mathcal{P}}
\newcommand{\ddd}  {\stackrel{\rm def}{=}}
\begin{document}
\title[Continuity of weighted operators \ldots]{Continuity of weighted operators, Muckenhoupt $A_p$ weights, and Steklov problem for orthogonal polynomials }
\author{Michel~Alexis, Alexander~Aptekarev, Sergey~Denisov}

\address{
\begin{flushleft}
Michel Alexis: malexis@wisc.edu\\\vspace{0.1cm}
University of Wisconsin--Madison\\  Department of Mathematics\\
480 Lincoln Dr., Madison, WI, 53706, USA\vspace{0.1cm}\\
\end{flushleft}
\vspace{0.3cm}
\begin{flushleft}
Alexander Aptekarev: aptekaa@keldysh.ru\\\vspace{0.1cm}
Keldysh Institute of Applied Mathematics\\ Russian Academy of Sciences\\
Miusskaya pl. 4, 125047 Moscow, RUSSIA\\
\end{flushleft}
\vspace{0.3cm}
\begin{flushleft}
Sergey Denisov: denissov@wisc.edu\\\vspace{0.1cm}
University of Wisconsin--Madison\\  Department of Mathematics\\
480 Lincoln Dr., Madison, WI, 53706, USA\vspace{0.1cm}\\
\vspace{0.1cm}
Keldysh Institute of Applied Mathematics\\ Russian Academy of Sciences\\
Miusskaya pl. 4, 125047 Moscow, RUSSIA\\
\end{flushleft}
\vspace{0.3cm}
}

\thanks{
The work of MA done in Sections 4 and 6 was supported by NSF grant
DMS-1147523.
The work of SD  is supported by
NSF-DMS-1464479,  NSF DMS-1764245, and Van Vleck Professorship Research
Award. SD gratefully
acknowledges the hospitality of IHES where part of this work was done.
}


\begin{abstract} We consider weighted operators acting on $L^p(\R^d)$ and show that they depend continuously on the weight $w\in A_p(\R^d)$ in the operator topology. Then, we use this result to estimate $L^p_w(\T)$ norm of polynomials orthogonal on the unit circle when the weight $w$ belongs to Muckenhoupt class $A_2(\T)$ and $p>2$. The asymptotics of the polynomial entropy is obtained as an application. 
\end{abstract}

\maketitle
\centerline{\it{To Peter Yuditskii on the occasion of his 65-th birthday}}
\setcounter{tocdepth}{1}

\tableofcontents



\section{Introduction}

Suppose $\mu$ is a probability measure on the unit circle $\T$ and $\{\phi_n(z,\mu)\}$ is the sequence of polynomials orthonormal with respect to $\mu$, i.e.
\begin{equation}\label{sad3}
\deg \phi_n=n,\qquad k_n\ddd {\rm coeff}_n\phi_n>0, \qquad  (\phi_n, \phi_k)_{L^2_\mu(\T)} = \delta_{n,k},
\end{equation}
where $\delta_{n,k}$ is the Kronecker symbol and ${\rm coeff}_jQ$ denotes the coefficient at the power $z^j$ in polynomial $Q$. One version of Steklov's problem in the theory of orthogonal polynomials  can be phrased as follows: given a  Banach space $X$ with norm $\|\cdot\|_X$, what regularity of $\mu$ is needed to have $\sup_{n\in \mathbb{N}}\|\phi_n(z,\mu)\|_X<\infty$? This problem has a long history. It goes back to Steklov's conjecture which asked to prove that the sequence $\{p_n(x,\rho)\}$ is bounded for every $x\in (a,b)$, where $\{p_n\}$ are polynomials orthonormal on the interval $[a,b]$ with respect to a weight $\rho$ that satisfies $\rho(x)\ge c>0, x\in [a,b]$. The negative answer to this question was given by Rakhmanov \cite{Rahmanov,Rahmanov2} and the sharp estimates on supremum norm were obtained only recently in \cite{adt}.  If $X=L^{2}_\mu(\T)$, we have $\|\phi_n\|_X=1$ by definition. In this paper, we will be concerned with the case  when  $X=L^p_\mu(\T), p>2$ and absolutely continuous $\mu$  is given by its weight, i.e., $d\mu=\frac{w}{2\pi} d\theta$. It is the natural choice since the space $L^p_w(\T)$ interpolates between the trivial case when $X=L^{2}_w(\T)$ and the space $L_w^\infty(\T)$, which was studied in \cite{adt,denik1} for  weights $w$ that satisfy Steklov's condition: $w^{-1}\in L^\infty(\T)$.

We recall the definition of Muckenhoupt class $A_p(\T)$ (see \cite{stein}, p.194).\smallskip

\noindent{\bf Definition.}
The weight $w\in A_p(\T), p\in (1,\infty)$ if 
\begin{equation}\label{sd_00}
[w]_{A_p(\T)}\ddd\sup_{I}\left(\langle w\rangle_I\left(\langle w^{\frac{1}{1-p}}\rangle_I\right)^{{p-1}}\right)<\infty,\,\quad \langle w\rangle_I\ddd\frac{1}{|I|}\int_Iwd\theta\,,
\end{equation}
where $I$ is an arc in $\T$. \smallskip

Given $w\in A_2(\T)$, we define the following quantity
\[
p_{\rm cr}(t)= \sup \{p: \sup_n\|\phi_n(z,w)\|_{L^p_w(\T)}<\infty, \,[w]_{A_2(\T)}\le t\}\,.
\]
Clearly, $p_{\rm cr}(t)$ is non-increasing on $[1,\infty)$ as a function in $t$  and $p_{\rm cr}(t)\ge 2$. The study of how $p_{\rm cr}(t)$ depends on $t$ amounts to considering another more precise version of Steklov's problem. Our first main result is the following theorem.
\begin{Thm}\label{t2}  We have
\[
p_{\rm cr}(t)>2, \quad \lim_{t\to 1}p_{\rm cr}(t)=+\infty, \quad \lim_{t\to\infty}p_{\rm cr}(t)=2\,.
\]
\end{Thm}
\noindent {\bf Remark.} In Appendix, we take $w$ as Fisher-Hartwig weight and prove $p_{\rm cr}(t)<C(t-1)^{-1/2}$ for $t\in (1,2]$. For $t>2$, the estimate $p_{\rm cr}(t)<2+Ct^{-1/6}$ will be obtained in the third section.\smallskip

The proof of this theorem in the perturbative regime, i.e., when $t$ is close to $1$, requires the following general result in the theory of weighted $L^p$ spaces. Consider spaces $L^p(\R^d)$ or $L^p(\T^d)$, $d\in \mathbb{N}$.  If $\mathcal{H}$ is a linear bounded operator from $L^p(\R^d)$ to itself, its operator norm will be denoted by $\|\mathcal{H}\|_{p,p}$. Suppose $w\in A_p(\R^d)$ and $H$ is a linear operator that satisfies weighted bound \begin{equation}\|w^{1/p}Hw^{-1/p}\|_{p,p}\le \cal{F}([w]_{A_p},p), \quad p\in (1,\infty) \label{sd_01}
\end{equation}
with some $p\in (1,\infty)$ and function $\cal{F}(t,p)$ which is  continuous in $t$ on $(1,\infty)$. In what follows, we do not need to know $\cal{F}$ explcitely. However,  $\cal{F}$ is known in many applications. For example, the Hunt-Muckenhoupt-Wheeden theorem (\cite{stein}, p.205) shows that $H$ can be taken as a singular integral operator and recent breakthrough on domination of singular integrals by sparse operators provides the sharp dependence of $\cal{F}$ on $[w]_{A_p}$. In particular, for a large class of singular integral operators, one can take
$
\cal{F}(t,p)=C(p)t^{\max(1,(p-1)^{-1})},
$
(see, e.g., \cite{NazLer}, p.264).

 Recall that $f\in {\rm BMO}(\R^d)$ if
\[
\|f\|_{{\rm BMO}(\R^d)}\ddd \sup_{B}\,\langle|f-\langle f\rangle_B|\rangle_B<\infty\,,
\]
where $B$ denotes a ball in $\R^d$ (see, e.g., p.140 in \cite{stein}).
The theorem that comes next is a slight improvement of a result by Pattakos and Volberg \cite{v1,v2}, see also the paper \cite{Pap_Pat} where the sublinear operators were treated.
\begin{Thm}\label{t1} Suppose $p\in (1,\infty)$, $[w]_{A_p(\R^d)}<\infty$, $\|f\|_{\rm BMO}<\infty$, and $H$ satisfies \eqref{sd_01}. Consider $w_\delta =we^{\delta f}$. Then, there is $\delta_0(p,[w]_{A_p},\|f\|_{\rm BMO})>0$ such that 
\[
\|w_\delta^{1/p}Hw_\delta^{-1/p}-w^{1/p}Hw^{-1/p}\|_{p,p}<|\delta|C(p,[w]_{A_p},\|f\|_{\rm BMO},\cal{F})
\]
for all $\delta: |\delta|<\delta_0$.
\end{Thm} 
 
 Two corollaries of theorem \ref{t2} are straightforward and we give their proofs in the end of section 3.  To state them, we need a few definitions. Given a weight $w$, define
\begin{equation}\label{rev_6}
q_{\rm cr}(w)= \sup \{q: \|w^{-1}\|_{L^q(\T)}<\infty\}\,.
\end{equation}
Clearly, if $w\in A_2(\T)$ then $q_{\rm cr}(w)>1$ and $\lim_{[w]_{A_2}\to 1}q_{\rm cr}(w)=\infty$ as follows from the definition of $A_p(\T)$ and inclusion of Muckenhoupt classes (see theorem 1 in \cite{Vasyunin} where the sharp bounds were obtained).\smallskip

\noindent {\bf Definition.} If $w\in L^1(\T)$ and it has finite logarithmic integral, i.e., $\log w\in L^1(\T)$, we define function $D$, the Szeg\H{o} function, as an outer function in $\D$ that satisfies
\begin{equation}\label{rev12}
|D|^2= w\,.
\end{equation}
The formula for $D$ is 
\begin{equation}\label{rev15}
D(z)=\exp\left(  \frac{1}{2\pi } \int_{\T} \frac{1+\bar\xi z}{1-\bar\xi z}\log \sqrt{ w(\theta)}d\theta\right),\, \xi=e^{i\theta}\,,z\in \D\,.
\end{equation}

\noindent {\bf Remark.} If $w\in A_2(\T)$, then $w^{-1}\in L^1(\T)$. Thus, $\log w\in L^1(\T)$ and $D$ is well-defined.\smallskip

 Given a polynomial $Q$ of degree at most $n$, its reversed polynomial  $Q^*$ is defined by $Q^*=z^n\overline{Q(1/\bar{z})}$. Notice that the map $Q\mapsto Q^*$ depends on $n$. Our first corollary establishes the asymptotics of $\{\phi_n^*\}$ (and thus of  $\{\phi_n\}$ since $\phi_n(\xi)=\xi^n\overline{\phi_n^*(\xi)}$ if $\xi\in \T$).

\begin{Cor}Suppose $[w]_{A_2}<\infty$ and $\|\frac{w}{2\pi}\|_1=1$, then 
\[
\lim_{n\to\infty}\|\phi_n^*-D^{-1}\|_{L^p_w(\T)}=0
\]
for every $p\in [2,\min \Bigl( p_{\rm cr}([w]_{A^2}),2(1+q_{\rm cr}(w))\Bigr))$.\label{cor_2}
\end{Cor}

Another application of theorem \ref{t2} has to do with the asymptotics of polynomial entropy $E(n,\mu)$, which is defined by
\[
E(n,\mu)=\int_\mathbb T |\phi_n(\xi,\mu)|^2\log|\phi_n(\xi,\mu)|d\mu\,,
\]
where $\xi=e^{i\theta},\, \theta\in [-\pi,\pi)$.

\begin{Cor}\label{entr}If $w\in A_2(\T)$, then
\[
\lim_{n\to\infty} E(n,w)=-\frac{1}{4\pi}\int_{-\pi}^\pi \log
w d\theta\,.
\]\label{cor_3}
\end{Cor}

Given a probability measure $\mu$ on $\T$, let $F$ be defined by
\begin{equation}\label{sd_51}
F(z)= \int_\T \frac{1+\bar \xi z}{1-\bar\xi z}d\mu,\quad \xi=e^{i\theta}\,.
\end{equation}
Notice that $\Re F>0$ in $\D$ and $F(0)=1$. For $\alpha\in \T$, consider the following one-parameter family (see, e.g., \cite{Simonbook}, p.36, formula (1.3.90))
\[
F_\alpha(z)\ddd \frac{\zeta + F(z)}{1+\zeta F(z)}\,,\quad  \zeta=\frac{1-\alpha}{1+\alpha}\in i(\R\cup \infty)\,.
\]
Function $F_\alpha$ also has positive real part in $\D$ and $F_\alpha(0)=1$, so
\[
F_{\rm \alpha}(z)= \int_\T \frac{1+\bar\xi z}{1-\bar\xi z}d\mu_{\rm \alpha}\,,
\]
which defines the family of Aleksandrov-Clark measures $\{\mu_{ \alpha}\}$. Taking $z=0$, we see that $\mu_\alpha$ is a  probability measure. If $\alpha=-1$, then $F_{-1}=1/F$ and the resulting measure is called dual for $\mu$, we will use notation $\mu_{\rm dual}(=\mu_{-1})$ for it. Measure $\mu_{\rm dual}$ plays an important role in the theory of polynomials orthogonal on the circle. In fact, the polynomials of the second kind $\{\psi_n\}$ defined by 
\[
\psi_n(z)=\int_\T \frac{1+z\bar\xi}{1-z\bar\xi}(\phi_n(\xi,\mu)-\phi_n(z,\mu))d\mu, \quad \xi=e^{i\theta}
\]
are orthonormal with respect to $\mu_{\rm dual}$ (see, e.g., \cite{Simonbook}, formulas (3.2.32) and (3.2.50) or section 1 in \cite{ger}).  The Muckenhoupt class $A_2(\T)$ turns out to be invariant with respect to taking dual. In fact, more general statement is true.
\begin{Thm}If $w\in A_2(\T)$ and $d\mu=\frac{w}{2\pi}d\theta$, then $\mu_{\alpha}$ is absolutely continuous and $d\mu_{\alpha}=\frac{w_{\alpha}}{2\pi}d\theta$ for every $\alpha\in \T$. Moreover, $w_{\alpha}\in A_2(\T)$.\label{t4}
\end{Thm}

This has an immediate implication for regularity of $\psi_n$. Indeed, if $w\in A_2(\T)$, then $d\mu_{\rm dual}=\frac{w_{\rm dual}}{2\pi}d\theta$ with $w_{\rm dual}\in A_2(\T)$ so theorem \ref{t2} can be applied and we get
\[
\sup_n\|\psi_n\|_{L^p_{w_{\rm dual}}(\T)}<\infty
\]
with  $p\in [2,p_{\rm cr}([w_{\rm dual}]_{A_2}))$. \smallskip

 The proofs of the main results in this paper involve complex interpolation, a suitable choice of the algebraic formulas, and a few facts from the general spectral theory.\smallskip
 
\noindent{\bf Previous results.} In \cite{adt}, it was proved  that, given every $q\in [1,\infty)$ and $n\in \mathbb{N}$, there is $w_\ast$ that satisfies $\|w_\ast\|_{L^q(\T)}<c_1, \|w_\ast^{-1}\|_{L^\infty(\T)}<c_2$ and nonetheless $\|\phi_n(\xi,w_\ast)\|_{L^\infty(\T)}\ge C(c_1,c_2,q) \sqrt n$ with parameters $c_1$ and $c_2$ being  $n$-independent. By Nikolskii inequality (see p.102, theorem 2.6, \cite{DeVore}), we see that $\|\phi_n(\xi,w_\ast)\|_{L^p(\T)}>C(c_1,c_2,p,q)n^{1/2-1/p}$ for every $p\in [2,\infty)$. Since the weight $w_\ast$ is bounded below by $c_2^{-1}$, one also gets $\|\phi_n(\xi,w_\ast)\|_{L^p_{w_\ast}(\T)}>C(c_1,c_2,p,q)n^{1/2-1/p}$. Therefore, the stated conditions on $w$, i.e., 
\[
\|w\|_{L^q(\T)}<c_1,\quad  \|w^{-1}\|_{L^\infty(\T)}<c_2,\quad q\in [1,\infty)
\]
do not provide the uniform in $n$ weighted $L^p$ estimates for polynomials if $p>2$ is fixed.
The question what regularity of $w$ is enough to have $\sup_{n}\|\phi_n\|_{L^p(\T)}<\infty$ or $\sup_{n}E(n,w)<\infty$ has been  addressed  in \cite{nevai,link1,link2,link3,denik1,denik2}.  The following theorem was proved in \cite{denik2}. \smallskip

\begin{Thm}[{\bf Denisov-Rush}, \cite{denik2}]\label{dr}
Let $s\dd \|w\|_{\rm BMO(\T)}<\infty$ and $t\dd \|w^{-1}\|_{\rm BMO(\T)}<\infty$. Then, there is $p(s,t)>2$ such that $\sup_{n}\|\phi_n(\xi,w)\|_{L^p(\T)}<\infty$. 
\end{Thm}
\noindent

 We will see later that theorem \ref{t2} implies theorem \ref{dr} and, in fact, gives a qualitatively stronger statement. It appears that $A_2$ regularity of $w$ is, to the best of our knowledge, the  weakest general condition that provides weighted $L^p$ estimates on $\{\phi_n\}$.
 
As far as theorem \ref{t1} is concerned, the continuity of operators in the weighted spaces with respect to a weight has been addressed previously. In \cite{v1,v2}, Pattakos and Volberg show that $A_\infty(\R^d)$ is a metric space with metric defined by
\[
d_*(w_1,w_2)\ddd \|\log w_1-\log w_2\|_{\rm BMO}\,.
\]
These two authors studied other properties of $A_\infty(\R^d)$ as a metric space and established, among other things, the Lipschitz continuity of $\|H\|_{L^p_w,L^p_w}$ in $w\in A_p(\R^d)$ for $H$ that satisfies \eqref{sd_01}.

\smallskip

\smallskip

The structure of our paper is as follows. The second section contains the proof of theorem \ref{t1} along with related information about the Muckenhoupt class. Theorem \ref{t2} and its corollaries are proved in the third section. The analysis of the Christoffel-Darboux kernel for the case when $w\in A_2(\T)$ is done in section four. In section five, we discuss Alexandrov-Clark measures and give proof of theorem~\ref{t4}.  The appendix contains an example of weight in the Fisher-Hartwig class for which the asymptotics of the polynomials is known. This provides an upper estimate for $p_{\rm cr}(t)$ in the regime when $t$ is close to $1$.

\subsection{Notation}

\begin{itemize}

\item If $p\in [1,\infty]$, the dual exponent is denoted by $p'=p/(p-1)$.

\item Given a set $A\subseteq \R^d$ (or $A\subseteq \T$), we will use notation $A^c$ for its complement, i.e., $A^c=\R^d\backslash A$ (or $A^c=\T\backslash A$).

\item Given  two Banach spaces $L^p(X,\mu)$, $L^q(Y,\nu)$, and a linear bounded operator $T:L^p(X,\mu)\rightarrow L^q(Y,\nu)$, its norm is denoted by $\|T\|_{p,q}$.

\item By $L^p _w(\T)$ we mean the space $L^p_\mu (\T)$ where $d\mu = w \frac{d \theta}{2 \pi}$.

\item If $f$ is locally integrable in $\R^d$ and $B$ is a ball, then 
\[
\langle f\rangle_B\ddd \frac{1}{|B|}\int_B fdx\,.
\]

\item Given function $f\in L^1(\T)$, we will write $\frak{h}(f)$ to denote the operator of harmonic conjugation \cite{Koos98}, i.e., 
\begin{equation}\label{rev10}
\frak{h}(f)=\widetilde f(\xi)=\lim_{r \to 1}\frac{1}{2\pi}\int_{\T}f(\zeta)Q_{r}(\zeta, \xi)\,d\theta, \quad 
Q_{r}(\zeta, \xi) = \Im \frac{1 + r\bar \zeta \xi}{1 - r\bar \zeta \xi}, \quad\zeta=e^{i\theta}, \quad \xi\in \T\,.
\end{equation}

\item Given a function $f\in L^1(\T)$, the Poisson integral is defined by (see \cite{Koos98}, pp.2--3)
\begin{equation}\label{rev111}
\P(f, z) = \frac{1}{2\pi}\int_{\T}\frac{1 - |z|^2}{|1 - \bar \zeta z|^2}f(\zeta)d\theta, \quad z \in \D, \quad \zeta=e^{i\theta}\,.
\end{equation}
The Cauchy integral over $\T$ is defined by (see \cite{Koos98}, p.35)
\begin{equation}\label{rev9}
\mathcal{C}(f,z)= \frac{1}{2\pi }\int_{\T}\frac{f(\zeta)}{1-\bar\zeta z}d\theta, \quad z \in \D, \quad \zeta=e^{i\theta}\,.
\end{equation}

\item For two non-negative functions
$f_{1}$ and $f_2$, we write $f_1\lesssim f_2$ if  there is an absolute
constant $C$ such that
\[
f_1\le Cf_2
\]
for all values of the arguments of $f_1$ and $f_2$. If the constant depends on a parameter $\alpha$, we will write $f_1\le_\alpha f_2$. We define $\gtrsim$
similarly and say that $f_1\sim f_2$ if $f_1\lesssim f_2$ and
$f_2\lesssim f_1$ simultaneously.

\item The symbol $C_c^\infty(\R^d)$ denotes the space of infinitely smooth function with compact support in $\R^d$.

\item Given two operators, $A$ and $B$,  we use the symbol $[A,B]=AB-BA$ for their commutator.
\end{itemize}

\section{Weighted operators are continuous in $w\in A_p(\R^d)$}

We start by recalling a  few basic facts from the theory of $A_p(\R^d)$ weights (see, e.g.,  \cite{korey} and \cite{stein}).  Given the definition \eqref{sd_00}, the limiting case when $p\to\infty$ leads to  $A_\infty(\R^d)$ which is characterized by (see, e.g., \cite{hru}) 
\begin{equation}\label{sd_1}
[w]_{A_\infty(\R^d)}\ddd \sup_B \left(\langle w\rangle_B \exp\Bigl( -\langle \log w\rangle_B\Bigr)\right)\,.
\end{equation}
The following results are well-known.
\begin{Lem}[see, e.g., \cite{stein}, p.218] If $\|f\|_{\rm BMO}<\infty$, then there is $\delta_1(\|f\|_{\rm BMO})>0$ such that
\[
[e^{\delta f}]_{A_\infty(\R^d)}\lesssim 1
\]
for all $\delta: |\delta|<\delta_1(\|f\|_{\rm BMO})$.\label{sd_7}
\end{Lem}
\begin{proof}From John-Nirenberg theorem (\cite{stein}, pp.145-146), we have
\begin{equation}
\sup_B\left(\langle e^{\delta|f-\langle f\rangle_B|}\rangle_B\right)\lesssim 1\label{sd_v3_1}
\end{equation}
provided $|\delta|<\delta_1(\|f\|_{\rm BMO})$. 
In \eqref{sd_1}, take $w=e^{\delta f}$, to get
\[
[w]_{A_\infty}= \sup_B \left(\langle 
e^{\delta (f-\langle f\rangle_B)}\rangle_B
\right)\lesssim 1
\]
by \eqref{sd_v3_1}.
\end{proof}
The proofs for the next two lemmas  are immediate corollaries from theorem $1_\infty$ and theorem 1 in \cite{Vasyunin}.

\begin{Lem}\label{sd_8} Suppose $w\in A_\infty(\R^d)$. For every $p\in (1,\infty)$, there is $\delta_2(p,[w]_{A_\infty(\R^d)})>0$ such that 
\begin{equation}\label{rev_1}
[w^\delta]_{A_p(\R^d)}<C(p,[w]_{A_\infty(\R^d)})
\end{equation}
for every $\delta: |\delta|<\delta_2(p,[w]_{A_\infty(\R^d)})$.
\end{Lem}

\noindent{\bf Remark.} The exact dependence of the right-hand side in \eqref{rev_1} on the parameters will not be needed in this paper so we are only using the symbol $C$.

\begin{Lem}Given $p\in (1,\infty)$ and $w\in A_p(\R^d)$, there is $\delta_3(p,[w]_{A_p(\R^d)})>0$ such that $[w^{1+\delta}]_{A_p(\R^d)}\le C(p,[w]_{A_p(\R^d)})$ for $\delta\in [0,\delta_3)$.\label{sd_2}
\end{Lem}

Given these lemmas, we claim that
\begin{Lem}\label{sd_10}
For every $p\in (1,\infty),f\in {\rm BMO}(\R^d),$ and $w\in A_p(\R^d)$, we have
\begin{equation}\label{sd_9}
[we^{\delta f}]_{A_p(\R^d)}\le C(p,[w]_{A_p(\R^d)},\|f\|_{\rm BMO})\,,
\end{equation}
if $\delta:|\delta|<\delta_4(p,[w]_{A_p(\R^d)},\|f\|_{\rm BMO})$.
\end{Lem}
\begin{proof}
Consider \eqref{sd_00}. Given $w$ and some nonnegative $w_0$ we use H\"older's inequality
\begin{eqnarray*}
\left(\int_{B}ww_0dx\right)\left(\int_{B}(ww_0)^{{1}/{(1-p)}}dx\right)^{p-1}\le\hspace{6cm}
\\
\left(\int_{B}w^\alpha dx\right)^{1/\alpha}\left(\int_{B} w_0^{\alpha'}dx \right)^{1/\alpha'}\left(\int_{B}w^{{\alpha}/{(1-p)}}dx\right)^{(p-1)/\alpha}\left(\int_{B}w_0^{{\alpha'}/{(1-p)}}dx\right)^{(p-1)/\alpha'}\,,
\end{eqnarray*}
where $\alpha'$ is dual to $\alpha$ and $\alpha>1$ is chosen such that $w^\alpha\in A_p(\R^d)$ (this choice is warranted by lemma~\ref{sd_2}). Now, if we let $w_0= e^{\delta f}$, then $w_0^{\alpha'}\in A_p(\R^d)$ for small $\delta$ thanks to lemma \ref{sd_7} and lemma \ref{sd_8}. This yields \eqref{sd_9}.
\end{proof}

\begin{Lem}If $p\in (1,\infty)$, $w\in A_p(\R^d)$, $f\in {\rm BMO}(\R^d)$, and $H$ satisfies \eqref{sd_01},   then
\begin{equation}\label{sad_el}
\|w^{1/p}[H,f]w^{-1/p}\|_{p,p}\le C(p,[w]_{A_p(\R^d)},\|f\|_{\rm BMO},\cal{F})
\end{equation}
and
\begin{equation}
\|w^{1/p}[f,[H,f]]w^{-1/p}\|_{p,p}\le C(p,[w]_{A_p(\R^d)},\|f\|_{\rm BMO},\cal{F}).\label{sad_ele}
\end{equation}
\label{sd_11}
\end{Lem}
\begin{proof}Given two test functions $u,v\in C_c^\infty(\R^d)$, define operator-valued function 
\[
G(z)\ddd w^{1/p}e^{zf}He^{-zf}w^{-1/p}
\]
and consider $\widehat G(z)=(G(z)u,v)$, where the inner product is in $L_2(\R^d)$. $\widehat G(z)$ is analytic in $z$ around the origin and we can
write Cauchy integral formula with $|z|<\epsilon$, when $\epsilon$ is small enough (and depends only on $p,[w]_{A_p(\R^d)},$ and $\|f\|_{\rm BMO}$):
\[
\widehat G(z)=\frac{1}{2\pi i}\int_{|\xi|=\epsilon}\frac{\widehat G(\xi)}{\xi-z}d\xi, \quad \widehat G'(0)=(w^{1/p}[f,H]w^{-1/p}u,v)=\frac{1}{2\pi i}\int_{|\xi|=\epsilon}\frac{\widehat G(\xi)}{\xi^2}d\xi, \]
so
\[ |(w^{1/p}[H,f]w^{-1/p}u,v)|\lesssim \epsilon^{-1}\max_{|\xi|=\epsilon}|\widehat G(\xi)|\,.
\]
For any point $z: |z|=\epsilon$ on the circle, we can apply lemma \ref{sd_10} and \eqref{sd_01} to choose $\epsilon(p,[w]_{A_p},\|f\|_{\rm BMO})$ such that $\max_{|\xi|=\epsilon}|\widehat G(\xi)|<C(p,[w]_{A_p},\|f\|_{\rm BMO},\cal{F})\|u\|_{p}\|v\|_{p'}$ (here $p'$ is dual to $p$). This implies \eqref{sad_el} by the standard duality argument, i.e., by employing an identity
\[
\|O\|_{p,p}=\sup_{u,v\in C^\infty_c(\R^d),\|u\|_p\le 1, |v\|_{p'}\le 1}|(Ou,v)|\,,
\]
which holds for every linear bounded operator $O$ and $p\in (1,\infty)$.

 The estimate \eqref{sad_ele} follows from \eqref{sad_el} by taking $H$ in \eqref{sad_el} as a commutator $[H,f]$ itself and using \eqref{sad_el}.
\end{proof}

{\it Proof of theorem \ref{t1}.} Consider analytic operator-valued function defined for $z:\Re z\in [0,1]$,
\[
F(z)=w^{1/p}\exp\left(\alpha z f/p\right)H \exp\left(-\alpha z f/p\right)   w^{-1/p}-w^{1/p}Hw^{-1/p}-z\frac{\alpha}{p}w^{1/p}[f,H]w^{-1/p}\,,
\]
where the parameter $\alpha$ will be chosen later, it will depend on $p,\|f\|_{\rm BMO}$, and $[w]_{A_p(\R^d)}$ only.
Consider rectangle $\Pi= \{z:|\Im z|<1, 0<\Re z<1\}$. 
We will estimate the operator norm of $F$ on $\partial \Pi$ as follows. If $z\in \{z:|\Im z|=1,\Re z\in [0,1]\}\cup \{z:\Re z=1, \Im z\in [-1,1]\}$, the estimate is straightforward:
\[
\|F(z)\|_{p,p}\le C(p,[w]_{A_p},\cal{F})+C(p,[we^{\alpha f}]_{A_p},\cal{F})\le C(p,[w]_{A_p},\|f\|_{\rm BMO},\cal{F}), \quad \alpha:|\alpha|<\alpha_4(p,[w]_{A_p},\|f\|_{\rm BMO})\,,
\]
where we first used  \eqref{sad_el}, \eqref{sd_01}, and then lemma \ref{sd_10}. Now, we take test functions $u,v\in C^\infty_c(\R^d)$ and consider $\widehat F(z)=(F(z)u,v)$.
It is anaytic in $\Pi$ and continuous on $\overline{\Pi}$. On the interval $z=i\xi, |\xi|<1$, have
\[
\widehat F(0)=0, \quad 
\widehat F'(0)=0,\quad
\partial_\xi \widehat F(i\xi)=\frac{i\alpha}{p}(w^{1/p}e^{i\alpha \xi f/p}[f,H]e^{-i\alpha\xi  f/p}w^{-1/p}u,v)-\frac{i\alpha}{p}w^{1/p}[f,H]w^{-1/p}\,,
\]
and 
\[
\partial^2_{\xi\xi} \widehat F(i\xi)=\left(\frac{i\alpha}{p}\right)^2(w^{1/p}e^{i\alpha \xi f/p}[f,[f,H]]e^{-i\alpha\xi  f/p}w^{-1/p}u,v)\,,
\]
\[
|\partial^2_{\xi\xi} \widehat F(i\xi)|\le C(p,[w]_{A_p},\|f\|_{\rm BMO},\cal{F}) \|u\|_p\|v\|_{p'} \quad 
\]
by lemma \ref{sd_11}. The Fundamental Theorem of Calculus gives
\[
\widehat F(i\xi)=\int_0^\xi \left(\int_0^\tau\partial^2_{\tau\tau} \widehat F(i\tau)d\tau \right)d\xi, \quad |\widehat F(i\xi)|\le \xi^2 C(p,[w]_{A_p},\|f\|_{\rm BMO},\cal{F})\|u\|_p\|v\|_{p'}.
\]
The last bound implies 
\[ \| F(i\xi)\|_{p,p}\le \xi^2 C(p,[w]_{A_p},\|f\|_{\rm BMO},\cal{F})\]
after we use duality argument. Notice that the function $|\widehat F|$
 is subharmonic in $\Pi$. Thus, by mean-value inequality, one has
\[
|\widehat F(\delta)|\le \left(\int_{\partial \Pi}|\widehat F(\xi)|d\omega_\delta(\xi)\right)\,,
\]
where $\omega_z(\xi)$ denotes the harmonic measure at point $z$ (see, e.g., \cite{gar_hm} p.13, formula (3.4)). By duality again, 
\[
\|F(\delta)\|_{p,p}\le \left(\int_{\partial \Pi}\|F(\xi)\|_{p,p}  d\omega_\delta(\xi) \right)\,.
\]
When $\delta\to 0$, measure $\omega_\delta(\xi)$ concentrates on the left side of $\partial \Pi$ around point $0$ and we have $\lim_{\delta\to 0}\|F(\delta)\|_{p,p}=0$. Putting the estimates together, we  can make it more precise. 
Recall that the harmonic measure on the upper half-plane $\C^+$ with the reference point $z$ is given by 
\[
\frac{1}{\pi}\frac{\Im z}{\Im^2 z+(\Re z-t)^2}\,, z\in \C^+, \, t\in \R\,.
\]
Consider a conformal map $\phi$ from  $\C^+$ to $\Pi$. For example, we can take $\phi$ as the following 
 Schwarz-Christoffel integral  (see \cite{ts}, p.181 and pp.188-189, formula (6-76)):
\[
\phi(z)=C\int_0^z \frac{d\eta}{\sqrt{(1-\eta^2)(1-k^2\eta^2)}}, \quad z\in \C^+,
\]
where $C$ and $k$ are constants that can be found explicitly and $k\in (0,1)$.  Under the inverse  map $\phi^{-1}$, the left side $\{i\xi, |\xi|<1\}$ of $\Pi$ goes to the interval $[-1,1]$ and its right side $\{1+i\xi, |\xi|<1\}$ goes to $[k^{-1},\infty)\cup (-\infty,-k^{-1}]$. Clearly, $\phi(0)=0$. Now, we obtain
\[
\int_{\partial \Pi}\|F(\xi)\|_{p,p}d\omega_\delta(\xi)    \lesssim \int_{\R}\frac{\delta}{\delta^2+t^2}\|F(\phi(t))\|_{p,p}dt\,
\]
where $\phi(t):\R\to \partial \Pi$.
Substituting the estimates for $\|F\|_{p,p}$ and using $|\phi(z)/z|\sim 1,\, |z|<0.5$, we get
\begin{eqnarray*}
\int_{\R}\frac{\delta}{\delta^2+t^2}\|F(\phi(t))\|_{p,p}dt \le \hspace{4cm}\\
C(p,[w]_{A_p},\|f\|_{\rm BMO},\cal{F})\left(\int_{-0.5}^{0.5}\frac{\delta t^2}{\delta^2+t^2} 
dt+ \int_{|t|>0.5}\frac{\delta}{\delta^2+t^2}dt\right)\le \\
C(p,[w]_{A_p},\|f\|_{\rm BMO},\cal{F})\delta\,.
\end{eqnarray*}
Finally, we get the statement of the theorem since
\[
w^{1/p}\exp\left(\alpha\delta f/p\right)H \exp\left(-\alpha\delta  f/p\right)   w^{-1/p}-w^{1/p}Hw^{-1/p}=F(\delta)+\delta\frac{\alpha}{p}w^{1/p}[f,H]w^{-1/p}   \,,
\]
and
\[\|F(\delta)\|_{p,p}\le C(p,[w]_{A_p},\|f\|_{\rm BMO},\cal{F})\delta\,,\] \[\|\frac{\alpha}{p}w^{1/p}[f,H]w^{-1/p}  \|_{p,p}\le C(p,[w]_{A_p},\|f\|_{\rm BMO},\cal{F})\,.\]
\qed\bigskip

\noindent {\bf Remark.} Clearly, the theorem holds if $A_p(\R^d)$ is replaced by $A_p(\T)$.\bigskip

\section{Steklov problem in the theory of orthogonal polynomials: $w\in A_2(\T)$ and bounds for $\|\phi_n(z,w)\|_{L^p_w(\T)}$}

This section contains the proofs of theorem \ref{t2} and its two corollaries. In the proof of theorem \ref{t2}, we will consider separately two cases: when $[w]_{A^2(\T)}\in [1,2)$  and when $[w]_{A^2(\T)}\ge 2$. It will be more convenient for us to work with monic orthogonal polynomials, which are defined as
\[
\Phi_n(z,\mu)=\frac{\phi_n(z,\mu)}{k_n}\,.
\]
If $w\in A_2(\T)$, then $w^{-1}\in L^1(\T)$ by definition. Thus, $\log w\in L^1(\T)$ as well. This means that $\mu: d\mu=\frac{w}{2\pi}d\theta$ belongs to Szeg\H{o} class of measures and, consequently, the sequence $\{k_n\}$ has a finite and positive limit (see \cite{ger}, section 2). More precisely, we have an estimate:
\begin{equation}\label{fact11}
\exp\left(\frac{1}{4\pi}\int_{\mathbb{T}} \log w
d\theta\right)\leq
\left|\frac{\Phi_n(z,w)}{\phi_n(z,w)}\right|\le 1, \quad
\forall z\in \mathbb{C},
\end{equation}
(see, e.g., \cite{denik1}). This bound shows that we can focus on estimating $\|\Phi_n(\xi,w)\|_{L^p_w(\T)}$. 

Later in the text, we will need to use the second resolvent identity which is contained in the following proposition.
\begin{Prop}
Suppose $X$ is an Banach space and $H,V$ are linear bounded operators from $X$ to $X$. Then,
\begin{eqnarray*}
(I+H+V)^{-1}=(I+H)^{-1}-(I+H+V)^{-1}V(I+H)^{-1}, \\
(I+H+V)^{-1}=(I+H)^{-1}(I+V(I+H)^{-1})^{-1}\,,
\end{eqnarray*}
provided the operators involved are well-defined and bounded in $X$. Moreover, assuming $\|V\|\cdot \|(I+H)^{-1}\|<1$, we get
\begin{equation}\label{sd_33}
\|(I+H+V)^{-1}\|\le  \frac{\|(I+H)^{-1}\|}{1-\|V\|\cdot \|(I+H)^{-1}\|}\,.
\end{equation}
Finally, if $\|V\|<1$, then 
\begin{equation}\label{ssa1}
\|(I+V)^{-1}\|\le \frac{1}{1-\|V\|}\,.
\end{equation}
\end{Prop}
\noindent The proof of this proposition is a straightforward calculation.
The following well-known lemma  (see, e.g., \cite{korey}, corollary 6) will be important later on.
\begin{Lem}\label{sd_ll}
If $[w]_{A_2(\T)}=1+\tau, \tau\in [0,1]$, then
\[
\|\log w\|_{\rm BMO}\lesssim \sqrt\tau\,.
\]
\end{Lem}

Let $\cal{P}_{n}$ denote the orthogonal $L^2(\T)$ projection to the frequencies $\{1,\ldots, e^{in\theta}\}$. Consider the perturbative regime, i.e., the case when $[w]_{A_2(\T)}=1+\tau$ and $\tau\in [0,1]$.

\begin{Lem} We have $\lim_{\tau \to 0} p_{\rm cr}(1+\tau)=\infty$.\label{l1}
\end{Lem}

\begin{proof} Fix any $p\ge  2$. We need to show that there is $\tau>0$ small enough so that $[\upsilon]_{A_2}<1+\tau$ implies 
\[
\sup_n \|\Phi_n(z,\upsilon)\|_{L^p_\upsilon(\T)}<\infty\,.
\]
 Our argument is based on a representation (see, e.g., \cite{denik2}, formula (8) for $\Phi_n^*$):
 \begin{equation}\label{sd_old1}
\Phi_n=z^n-\upsilon^{-1}[\cal{P}_{n-1}, \upsilon]\Phi_n\,.
\end{equation}
This formula can be obtained by combining trivial identity $\Phi_n=z^n+\cal{P}_{n-1}\Phi_n$, which holds for all monic polynomials of degree $n$, with $\cal{P}_{n-1}(\upsilon \Phi_n)=0$, which follows from that fact that $\Phi_n$ is orthogonal to $\{1,z,\ldots,z^{n-1}\}$ in $L^2_\upsilon(\T)$.
Thus, we infer from \eqref{sd_old1} that
\[
\Bigl(\upsilon^{1/p} \Phi_n\Bigr)=\upsilon^{1/p}
z^n-\upsilon^{-1/p'}\cal{P}_{n-1}\upsilon^{1/p'}
\Bigl(\upsilon^{1/p}\Phi_n\Bigr)+\upsilon^{1/p}\cal{P}_{n-1}\upsilon^{-1/p}\Bigl(\upsilon^{1/p}\Phi_n\Bigr)\,.
\]
Denoting $\zeta_n\ddd \upsilon^{1/p}\Phi_n$, $O_{1,n}\ddd \upsilon^{-1/p'}\cal{P}_{n-1}\upsilon^{1/p'}-\cal{P}_{n-1}$,
$O_{2,n}\ddd \upsilon^{1/p}\cal{P}_{n-1}\upsilon^{-1/p}-\cal{P}_{n-1}$, we rewrite it as
\begin{equation}\label{sd_22}
\zeta_n=\upsilon^{1/p}z^n-O_{1,n}\zeta_n+O_{2,n}\zeta_n\,.
\end{equation}
If $\cal{P}^+$ denotes the orthogonal $L^2(\mathbb{T})$ projection   onto Hardy space $H^2(\mathbb{T})$ (Riesz projection), then we can write an identity
\begin{equation}\label{sd_13}
\cal{P}_{n}=\cal{P}^+-z^{n+1} \cal{P}^+ z^{-(n+1)}=z^{n+1}[z^{-(n+1)},\cal{P}^+]\,.
\end{equation}
We now apply theorem \ref{t1} with $H=\cal{P}^+$, $w=1$, and $w_\delta=e^{\delta f}=\upsilon$. Then, $f=\delta^{-1}\log \upsilon$ and lemma \ref{sd_ll} gives
\[
\|f\|_{\rm BMO}\lesssim \delta^{-1}\sqrt\tau\le 1\,,
\]
when $\tau<\delta^2$. Since $\|w^{1/p}\cal{P}^+w^{-1/p}\|_{p,p}\le \cal{F}([w]_{A^p},p)$ by Hunt-Muckenhoupt-Wheeden theorem,  the theorem~\ref{t1} then yields
\[
\lim_{\tau\to 0}\|\upsilon^{1/p}\cal{P}^+\upsilon^{-1/p}-\cal{P}^+\|_{p,p}=0
\]
for every $p\in (1,\infty)$. In particular, it also holds for $p'$: 
\[
\lim_{\tau\to 0}\|\upsilon^{1/p'}\cal{P}^+\upsilon^{-1/p'}-\cal{P}^+\|_{p',p'}=0\,.
\]
Indeed, we use the standard identity in the operator theory, which follows from duality considerations:
\[
\|\cal{O}\|_{p,p}=\|\cal{O}^*\|_{p',p'}\,,
\]
where $\cal{O}^*$ is adjoint operator to $\cal{O}$ with respect to $L^2$ inner product and $\cal{O}$ is linear bounded operator in $L^p$ space. Since $\cal{P}^+$ is self-adjoint in $L^2(\T)$, 
 we get \[\|\upsilon^{1/p'}\cal{P}^+\upsilon^{-1/p'}-\cal{P}^+\|_{p',p'}=\|\upsilon^{-1/p'}\cal{P}^+\upsilon^{1/p'}-\cal{P}^+\|_{p,p}\]
and hence
\[
\lim_{\tau\to 0}\|\upsilon^{-1/p'}\cal{P}^+\upsilon^{1/p'}-\cal{P}^+\|_{p,p}=0\,.
\]
Summarizing, \eqref{sd_13} gives two bounds
\[
\|O_{1,n}\|_{p,p}\le 2\|\upsilon^{-1/p'}\cal{P}^+\upsilon^{1/p'}-\cal{P}^+\|_{p,p},\, \|O_{2,n}\|_{p,p}\le 2\| \upsilon^{1/p}\cal{P}^+\upsilon^{-1/p}-\cal{P}^+   \|_{p,p}
\]
that hold uniformly in $n$. Therefore, 
\[
\lim_{\tau\to 0}\|O_{2,n}\|_{p,p}=0, \quad \lim_{\tau\to 0}\|O_{1,n}\|_{p,p}=0\,.
\]
 Now, we  apply \eqref{ssa1}  with $V=O_{1,n}$ to \eqref{sd_22} in the space $L^p(\T)$. This gives the statement of the lemma. Here, we notice that $\sup_{n}\|z^n\upsilon^{1/p}\|_p<\infty$ because $\upsilon \in A_2(\T)\subset L_1(\T)$. 
\end{proof}

Next, we consider more complicated case when $[w]_{A_2(\T)}\ge 2$.

{\bf Remark.} We have $w^{-1/p'}=(w^{-p/p'})^{1/p}$ and
\begin{equation}\label{sad1}
[w^{-p/p'}]_{A_{p}(\T)}=[w]_{A_{p'}(\T)}^{p/p'}
\end{equation}
as can be directly verified.\smallskip

\begin{Lem}\label{rev_3} For every $w\in A_2(\T)$ and $l\in \mathbb{N}$,  define  a simple function $w_l$ as follows: let $w_l=\langle w\rangle_{I_j}$ on each  interval $I_j=2^{-l}(2\pi)[j,j+1), j=0,\ldots,2^l-1$. Then, $\lim_{l\to\infty}\Phi_n(z,w_l)=\Phi_n(z,w)$ uniformly in $z$ over compacts in $\C$ and  \[[w_l]_{A_2(\T)}\le C([w]_{A_2(\T)}).\]
\end{Lem}
\begin{proof}
From the construction, we immediately get
$\{w_l\}{\to} w$ in the weak--$(\ast)$ sense when $l\to\infty$. Since the coeffcients of $\Phi_n(z,\mu)$ depend continuously on the moments of measure $\mu$, we have the first statement of the lemma. The second one can be verified directly using the definition of $A_2(\T)$ characteristic.

\end{proof}

Next, we need the following interpolation result. Given $w\in A_2(\T)$ and $p_*\geq 2$, define
\begin{equation}\label{sa1}
Q_{w,p(z)}\ddd w^{-1/p'(z)}\cal{P}_{n-1}w^{1/p'(z)}-w^{1/p(z)}\cal{P}_{n-1}w^{-1/p(z)}\,,
\end{equation}
where
\begin{equation}\label{sad7}
\frac{1}{p(z)}=\frac{z}{p_*}+\frac{1-z}{2}, \quad \frac{1}{p'(z)}=1-\frac{1}{p(z)}=\frac{1+z}{2}-\frac{z}{p_*}, \quad \Re z\in [0,1]\,,
\end{equation}
so that $1/p(z)+1/p'(z)=1$.

\begin{Prop}\label{ss1} Suppose $w,w^{-1}\in L^\infty(\T)$, parameter $\kappa$ is real, and 
\begin{equation}\label{rev_2}
 \sup_{0\le \Re z\le 1}\|Q_{w,p(z)}\|_{p(t),p(t)}<\infty\,,
 \end{equation}
where $t\ddd \Re z\in [0,1]$. If there is a positive number ${\bf \Lambda}$ such that
\[
\|(I-\kappa Q_{w,p(t+iy)})^{-1}\|_{p(t),p(t)}\le 2{\bf \Lambda}
\]
for all $t\in [0,1]$ and $y\in \R$, then there is an $t_*({\bf \Lambda})\in (0,1]$,  so that
\[
\|(I-\kappa Q_{w,p(t+iy)})^{-1}\|_{p(t),p(t)}\le {\bf \Lambda}
\]
for all $y\in \R$ and $t\in [0,t_*]$.
\end{Prop}

\begin{proof}
We notice that $Q_{w,p(iy)}$ is bounded and antisymmetric operator in Hilbert space $L^2(\T)$. Therefore, $\|(I-\kappa Q_{w,p(iy)})^{-1}\|_{2,2}\le 1$.  Given conditions $w,w^{-1}\in L^\infty(\T)$, it is easy to check that the operator-valued function  $(I-\kappa Q_{w,p(z)})^{-1}$ is analytic  and continuous in the sense of Stein (p.209, \cite{bsh}). Applying Stein's interpolation theorem, we get
\[
\|(I-\kappa Q_{w,p(t+iy)})^{-1}\|_{p(t),p(t)}\le \exp\left(
\frac{\sin(\pi t)}{2}\int_\R \frac{\log (2{\bf \Lambda)}}{\cosh (\pi y)+\cos(\pi t)}dy\right)=1+O(t), \quad t\to 0\,,
\]
which proves the proposition.
\end{proof}
\noindent{\bf Remark.} We emphasize here that positive $t_*$ does not depend on $n$ or $w$.

Now, we are ready to prove the following lemma.

\begin{Lem} For every $t\ge  2$, we have $p_{\rm cr}(t)>2$. \label{l2}
\end{Lem}

\begin{proof} 
Consider $w\in A_2(\T)$.
It will be more convenient later on to work with weights which are bounded above and below. With fixed $n$, we can use lemma \ref{rev_3} to approximate $w$ by $w_n$ which satisfies
\[
\|w_n\|_{L^\infty(\T)}<C(n,w),\quad \|w_n^{-1}\|_{L^\infty(\T)}<C(n,w)\,,
\]
\[
[w_n]_{A_2(\T)}\le \gamma\ddd C([w]_{A_2}), \, n\in \mathbb{N}
\]
and
\[
 \quad |\Phi_n(z,w)|\le 2|\Phi_n(z,w_n)|
\]
for each $z\in \T$. In what follows, we suppress the dependence of $w_n$ in $n$ and  do the proof understanding that $w$ depends on $n$ and satisfies
\[
\|w\|_{L^\infty(\T)}<\infty,\quad \|w^{-1}\|_{L^\infty(\T)}<\infty, \quad [w]_{A_2(\T)}\le \gamma<\infty\,,
\]
where $\gamma$ does not depend on $n$.

As in the proof of  lemma \ref{l1}, we can write 
\[
\zeta_n=w^{1/p}z^n+Q_{w,p}\zeta_n\,,
\]
where $\zeta_n \ddd w^{1/p}\Phi_n$ and $Q_{w,p}\ddd -B_n+C_n, B_n\ddd w^{-1/p'}\cal{P}_{n-1}w^{1/p'}, C_n\ddd w^{1/p}\cal{P}_{n-1}w^{-1/p}$ and all operators are considered in Banach space  $L^p(\T)$.  It is sufficient to prove that 
\begin{equation}\label{sd_m}
 \sup_n\|(I-Q_{w,\widetilde p_\gamma})^{-1}\|_{\widetilde p_\gamma,\widetilde p_\gamma}<\infty
\end{equation}
with some $\widetilde p_\gamma>2$ because $\sup_n\|w^{1/p}z^n\|_{ p}<\infty$ and 
\[
\zeta_n=(I-Q_{w,p})^{-1}(w^{1/p}z^n)\,.
\]
 By open inclusion of Muckenhoupt classes (see \cite{stein}, corollary on p.202 or theorem 1 in \cite{Vasyunin}), there is $\widehat p_\gamma>2$ such that 
${\widehat p}'_\gamma<2$ and $
\widehat\gamma\ddd [w]_{A_{\widehat p'_\gamma}}<\infty\,.
$ Thus, by \eqref{sad1}, 
\begin{equation}\label{sad11}
[w^{-p/p'}]_{A_{p}}=[w]_{A_{p'}}^{p/p'}\le \widehat \gamma^{\widehat p/\widehat p'}
\end{equation}
for all $p\in [2,\widehat p_\gamma]$. We need this bound to  control $B_n$ through writing it as 
\[
B_n=(w^{-p/p'})^{1/p}\cal{P}_{n-1}(w^{-p/p'})^{-1/p}
\]
and viewing $w_1\ddd w^{-p/p'}$ as element of $A_p(\T)$. Now, we use Hunt-Muckenhoupt-Wheeden theorem, which implies that 
\begin{equation}\label{sd_012}
\sup_n \|B_n\|_{p,p}=\sup_n\|w_1^{1/p}\cal{P}_{n-1}w_1^{-1/p}\|_{p,p}<\cal{F}_1(p,\gamma)\,,
\end{equation}
where $\cal{F}_1$ is defined for $p\in [2,\widehat p_\gamma]$. Analogous bound for $C_n$ is obvious:
\begin{equation}\label{sd_022e}
\sup_n\|C_n\|_{p,p}<\cal{F}_2(p,\gamma)
\end{equation}
for all $p\in (2,\infty)$ since $w\in A_2(\T)\subset A_p(\T)$.
Define $Q_{w,p(z)}$  by \eqref{sa1} and take $p_*\in [2,\widehat p_\gamma]$.
The bounds \eqref{sd_012} and \eqref{sd_022e} imply that  
\[
\sup_n \|Q_{w,p(z)}\|_{p(t),p(t)}<\infty 
\]
for $t=\Re z\in [0,1]$. 

Now, we proceed as follows. Recall, see \eqref{sd_m}, that our goal is to show  that  
$
(I-Q_{w,\widetilde p_\gamma})^{-1}
$ is bounded in $L^{\widetilde p}(\T)$ for some $\widetilde p_\gamma>2$ with bound on the operator norm independent in $n$. In \eqref{sa1}, we take parameter $p_*$ as follows: $p_*^{(1)}=\widehat p_\gamma$ and define $p_1(z)\ddd p(z)$ where $p(z)$ is from  \eqref{sad7}. Consider $Q^{(j)}_{w,p(z)}\ddd jQ_{w,p(z)}/N, j=1,\ldots,N$ where $N$ is large and will be fixed later (it will depend on $\gamma$ only). Notice that, by \eqref{sd_012} and \eqref{sd_022e}, we get
\[
\|Q_{w,p(t+iy)}\|_{p(t),p(t)}\le \|w^{-1/p'(t)}\cal{P}_{n-1}w^{1/p'(t)}\|_{p(t),p(t)}+\|w^{1/p(t)}\cal{P}_{n-1}w^{-1/p(t)}\|_{p(t),p(t)}<C_\gamma\,.
\]
Let ${\bf \Lambda}$ be an absolute constant larger than one.
We take $N$ to satisfy 
\begin{equation}\label{sd_44}
 1-C_\gamma{\bf \Lambda}/N>1/2\,.
\end{equation}
Next, we use \eqref{ssa1} to get
\[
\|(I-Q^{(1)}_{w,p(t+iy)})^{-1}\|_{p(t),p(t)}\le \frac{1}{1-C_\gamma /N}\le \frac{1}{1-C_\gamma {\bf \Lambda}/N}\le 2\le  2{\bf \Lambda}
\]
since ${\bf \Lambda}>1$ by our choice. We continue with an inductive argument in which the bound for  $\{Q^{(j)}_{w,p(z)}\}$ provides the bound for $\{Q^{(j+1)}_{w,p(z)}\}$ when  $j=1,\ldots,N-1$. \smallskip

{$\bullet$ \bf Base of induction: handling  $Q^{(1)}_{w,p(z)}$.} Apply proposition \ref{ss1} with $\kappa=1/N$ to get an absolute constant $t_*$ so that
\[
\|(I-Q^{(1)}_{w,p(t+iy)})^{-1}\|_{p(t),p(t)}\le {\bf \Lambda}
\]
for $t\in [0,t_*]$ and $y\in \R$. Next, we use \eqref{sd_33} with $H=-Q^{(1)}_{w,p(t+iy)}$ and $V=-N^{-1}Q_{w,p(t+iy)}$. This gives
\begin{equation}\label{sa4}
\|(I-Q^{(2)}_{w,p(t+iy)})^{-1}\|_{p(t),p(t)}\le \frac{{\bf \Lambda}}{1-C_\gamma {\bf \Lambda}/N}\le 2{\bf \Lambda}, \quad t\in [0,t_*]
\end{equation}
by \eqref{sd_44}.

That finishes the first step. Next, we will explain how estimates on $Q^{(2)}_{w,p(z)}$ give bounds for $Q^{(3)}_{w,p(z)}$.

{$\bullet$ \bf Handling  $Q^{(2)}_{w,p(z)}$.} In proposition \ref{ss1}, we now take  $\kappa=\kappa_2\ddd 2/N, p^{(2)}_*\ddd p_1(t_*)=p(t_*)$ (here $p(t_*)$ is obtained at the previous step) and compute new $p_2(z),p_2'(z)$ by \eqref{sad7}:
\begin{equation}\label{sad66}
\frac{1}{p_2(z)}=\frac{z}{p(t_*)}+\frac{1-z}{2}=\frac{zt_*}{p_*}+\frac{1-zt_*}{2}=\frac{1}{p_1(zt_*)}=\frac{1}{p(zt_*)}\,.
\end{equation}
Therefore, when $z$ belongs to $0<\Re z<1$, $zt^*$ belongs to $0<\Re z<t_*$ and $p_2(z)=p(zt_*)$. In this domain, we have an estimate \eqref{sa4} which can be rewritten as
\[
\|(I- Q^{(2)}_{w,p_2(t+iy)})^{-1}\|_{p_2(t),p_2(t)}\le 2{\bf \Lambda}, \quad t\in [0,1],\quad y\in \R\,,
\]
where $p_2(z)$ is different from $p_1(z)=p(z)$ only by the choice of parameter  $p_*$ in \eqref{sad7} and is in fact a rescaling of the original $p(z)$ as follows from  \eqref{sad66}. Thus, from proposition \ref{ss1}, we have 
\[
\|(I-Q^{(2)}_{w,p_2(t+iy)})^{-1}\|_{p_2(t),p_2(t)}\le {\bf \Lambda}
\]
for $t\in [0,t_*],y\in \R$. We use the perturbative bound \eqref{sd_33} one more time with $H=-Q^{(2)}_{w,p_2(t+iy)}$ and $V=-N^{-1}Q_{w,p_2(t+iy)}$ to get
\[
\|(I-Q^{(3)}_{w,p_2(t+iy)})^{-1}\|_{p_2(t),p_2(t)}\le 2{\bf \Lambda}
\]
for $t\in [0,t_*],y\in \R$.

{$\bullet$ \bf Induction in $j$ and  the bound for $Q^{(N)}_{w,p(z)}$.} Next, we take $p_*^{(3)}\ddd p_*^{(2)}(t_*)$ and repeat the process in which the bound 
\[
\|(I- Q^{(j)}_{w,p_j(t+iy)})^{-1}\|_{p_j(t),p_j(t)}\le 2{\bf \Lambda}, \quad t\in [0,1],\quad y\in \R\,,
\]
implies 
\[
\|(I-Q^{(j+1)}_{w,p_{j+1}(t+iy)})^{-1}\|_{p_{j+1}(t),p_{j+1}(t)}\le 2{\bf \Lambda}
\]
for $t\in [0,1]$ and $y\in \R$.  Notice that each time the new $p_j(z)$ is in fact a rescaling of the original $p(z)$ by $t_*^{j-1}$ as can be seen from a calculation analogous to  \eqref{sad66}. In $N-1$ steps, we get
\[
\|(I-Q^{(N)}_{w,p_{N-1}(t+iy)})^{-1}\|_{p_{N-1}(t),p_{N-1}(t)}\le  2{\bf \Lambda}\,, \quad t\in [0,t_*], \quad y\in \R.
\]
Thus, taking $y=0$ and $t=t_*$, and recalling that $p_{N-1}(z)=p(t_*^{N-2}z)$, one has
\[
\|(I-Q^{(N)}_{w,p(t_*^{N-1})})^{-1}\|_{p(t_*^{N-1}),p(t_*^{N-1})}\le  2{\bf \Lambda}\,.
\]
Since $Q^{(N)}_{w,p(t_*^N)}=Q_{w,p(t_*^N)}$, we get \eqref{sd_m} with 
\[
\widetilde p_\gamma=\frac{2\widehat p_\gamma}{2t_*^{N-1}+\widehat p_\gamma(1-t_*^{N-1})}.
\] 
The estimates \eqref{sd_44} implies that we can take
$
N\sim C_\gamma\,.
$
\end{proof}
{\it Proof of theorem \ref{t2}.}
From lemma \ref{l1} and lemma \ref{l2}, we get that $p_{\rm cr}(t)>2$ and $\lim_{t\to 1} p_{\rm cr}(t)=\infty$. To show that $p_{\rm cr}(t)\to 2$ when $t\to\infty$, it is enough to start with  arbitrarily large $t$ and present a weight $\widehat w$  such that $[\widehat w]_{A_2(\T)}\le t$ and $\sup_n \|\phi_n(\xi,\widehat w)\|_{L^{p(t)}_{\widehat w}(\T)}=+\infty$ with some $p(t)$ which depends on $t$ and $\lim_{t\to\infty}p(t)=2$. To this end, we use the following result established in \cite{denik1}, theorem~3.2: {\it given any $t>2$, there is a weight $w$ that satisfies $1\le w\le t$ and a subsequence $\{k_n\}$ such that
\[
 \quad \|\phi_{k_n}(\xi,w)\|_{L^\infty(\T)}\ge C(t)k_n^{1/2-ct^{-1/6}}\,.
\]
}
The weight $w$ in the statement does not  satisfy condition $\|\frac{w}{2\pi}\|_{L^1(\T)}=1$. However, for $\widehat w=2\pi w/\|w\|_{L^1(\T)}$, we will have 
\begin{equation}
\Bigl\|\frac{\widehat w}{2\pi}\Bigr\|_{L^1(\T)}=1,\quad
\frac{\sup_\T \widehat w}{\inf_\T \widehat w}\le t
\end{equation}
and
\[
\|\phi_{k_n}(\xi,\widehat w)\|_{L^\infty(\T)}\ge C(t)k_n^{1/2-ct^{-1/6}}\,.
\]
Nikolskii inequality (see p.102, theorem 2.6, \cite{DeVore}) gives
$
\|\phi_{k_n}(\xi,\widehat w)\|_{L^p(\T)}\ge C(t,p)k_n^{1/2-1/p-ct^{-1/6}}
$
and thus
\[
\|\phi_{k_n}(\xi,\widehat w)\|_{L^p_{\widehat w}(\T)}\ge C(t,p)k_n^{1/2-1/p-ct^{-1/6}}\,.
\]
The weight $\widehat w$ satisfies the trivial bound
$
[\widehat w]_{A_2(\T)}\le t
$. Therefore, 
\[
p_{\rm cr}(t)\le \frac{2t^{1/6}}{t^{1/6}-2c}=2+O(t^{-1/6}), \, t\to\infty \,.
\]
\qed\smallskip

\noindent{\bf Remark.} Some lower bounds on $p_{\rm cr}(t)$ when $t\to 1$ and $t\to\infty$ can be traced through the proof.   We do not include these calculations here.  \smallskip

{\it Proof of corollary \ref{cor_2}.}
We have (see \cite{KH01}, formula (5.37) or \cite{ger}, section 2)
\begin{equation}\label{sd_56}
\lim_{n\to\infty}\|\phi_n^*-D^{-1}\|_{L^2_w(\T)}=0\,.
\end{equation}
Recall that $q_{\rm cr}(w)$ was defined in \eqref{rev_6}.
Take $\widetilde p\in [2,\min (p_{\rm cr}([w]_{A^2}),2(1+q_{\rm cr}(w))))
$. For $p\in [2,\widetilde p)$,
  we use H\"older's inequality
\begin{equation}\label{sd_98}
\int_{\T}|\phi_n^*-D^{-1}|^pw d\theta\le \left(\int_\T |\phi_n^*-D^{-1}|^{p_1\alpha}wd\theta\right)^{1/\alpha}\cdot \left(\int_\T |\phi_n^*-D^{-1}|^{p_2\alpha'}wd\theta\right)^{1/\alpha'}\,,
\end{equation}
where $p_1+p_2=p, p_1\alpha=\widetilde p, p_2\alpha'=2, \alpha^{-1}+\alpha'^{-1}=1, \alpha\in (1,\infty)$. In fact, 
solving these equations gives $\alpha=(\widetilde p-2)/(p-2)$, $p_1=\widetilde p(p-2)/(\widetilde p-2)$, $p_2=2(\widetilde p-p)/(\widetilde p-2)$. The second factor in the right hand side of \eqref{sd_98} converges to zero due to \eqref{sd_56}. For the first one, we apply the triangle inequality to write
\[
\sup_n\left(\int_\T |\phi_n^*-D^{-1}|^{\widetilde p}wd\theta\right)^{1/\widetilde p}\le \sup_{n}\|\phi_n^*\|_{\widetilde p,w}+\|D^{-1}\|_{\widetilde p,w}\,.
\]
The first term is finite thanks to theorem \ref{t2}. For the second one, we use $w=|D|^2$ to write
\[
\|D^{-1}\|_{\widetilde p,w}^{\widetilde p}=\int_\T |D^{-1}|^{\widetilde p}wd\theta=\int_\T w^{1-\widetilde p/2}d\theta<\infty\,,
\]
because $\widetilde p/2-1<q_{\rm cr}(w)$.
\qed\bigskip

{\it Proof of corollary \ref{cor_3}.} Let $S\ddd D^{-1}$ for shorthand.
Recall that $|\phi_n|=|\phi_n^*|$ on $\T$.
The following inequality follows from the Mean Value Formula
\[
|x^2\log x-y^2\log y|\lesssim (1+x|\log x|+y|\log y|)|x-y|, \quad x,y\ge
0\,.
\]
Hence,
\[
\int_{-\pi}^\pi
||\phi^*_n|^2\log|\phi^*_n|-|S|^2\log|S||wd\theta\lesssim
\int_{-\pi}^\pi(1+|\phi^*_n\log|\phi^*_n||+|S\log|
S||)||\phi_n^*|-|S||wd\theta\,.
\]
Then, one can write
\begin{eqnarray*}
\int_{-\pi}^\pi(1+|\phi^*_n\log|\phi_n||+|S\log|
S||)||\phi_n^*|-|S||wd\theta\le \hspace{4cm}\\ C(\delta)  \left(\int_{-\pi}^\pi(1+|\phi^*_n|^{2+\delta}+|S|^{2+\delta})wd\theta\right)^{1/2} \left(\int_{\T}|\phi_n^*-S|^2wd\theta\right)^{1/2}
\end{eqnarray*}
by applying Cauchy-Schwarz inequality and the trivial bound: $(1+u|\log u|)^2\le C(\delta)(1+u^{2+\delta}),\,\, \delta>0$. The second factor converges to zero when $n\to\infty$ due to \eqref{sd_56}. For the first one, theorem \ref{t2} and identity $|S|=w^{-1/2}$ allow us to find $\delta>0$ such that 
\[
\sup_{n}\int_{-\pi}^\pi(|\phi^*_n|^{2+\delta}+|S|^{2+\delta})wd\theta<\infty\,.
\]
\qed

In the rest of this section, we will show that theorem \ref{t2} implies theorem \ref{dr}. We start with the following lemma.
\begin{Lem} If $w, w^{-1}\in {\rm BMO}(\T)$, then $w\in A_2(\T)$.
\end{Lem}
\begin{proof}
Let $s\ddd \|w\|_{{\rm BMO}(\T)}, t\ddd \|w^{-1}\|_{{\rm BMO}(\T)}$ for shorthand. 
Consider any interval $I\subseteq \T$. We define $a\ddd \langle w\rangle_I, b\ddd \langle w^{-1}\rangle_I$.  We have
\[
\langle |w-a|\rangle_I\le s, \quad\langle|w^{-1}-b|\rangle_I\le t
\]
by the definition of BMO space.
To estimate $A_2(\T)$ characteristic, we need to bound $ab$. We assume without loss of generality that $I=[0,1]$ and that $a\le b$.  Apply triangle's inequality and an estimate 
\[
\frac{1}{|I|}\|w-\langle w\rangle_I\|_{L^2(I)}^2\lesssim s^2
\]
(see \cite{stein}, p.144, formula (7)), to get
 \begin{equation}
 \|w\|_2\le \|w-a\|_2+\|a\|_2\lesssim s+a\,,
 \end{equation}
where here and in the rest of the proof all estimates are done with respect to $I=[0,1]$.
Consider a set $\Omega\ddd \{|w^{-1}-b|\le 0.5b\}$. By John-Nirenberg inequality (\cite{stein}, p.145, formula (8)), we can estimate the measure of its complement via
\begin{equation}\label{sd_kek}
 |\Omega^c|\lesssim \exp \left(-c_1b t^{-1}\right)\,,
\end{equation}
where $c_1$ is an absolute positive constant.
We can rewrite $\Omega$ as follows $\Omega=\{0.5b\le w^{-1}\le 1.5b\}=\{2/(3b)\le w\le 2/b\}$ and this formula shows that
\begin{equation}\label{sade1}
\int_{w>2/b}d\theta\le |\Omega^c| \lesssim \exp(-c_1bt^{-1})\,.
\end{equation}
Then,
\[
a=\int_{w\le 2/b}wd\theta+\int_{w>2/b}wd\theta
\]
and consequently
\[
\int_{w>2/b   }wd\theta=a-\int_{w \le 2/b }wd\theta\ge a-2/b\,.
\]
On the other hand, by Cauchy-Schwarz inequality and \eqref{sade1},
\[
\int_{w>2/b   }wd\theta \le \|w\|_2 \left(\int_{w>2/b}d\theta\right)^{1/2}\lesssim 
(s+a)\exp \left(-c_1bt^{-1}/2\right)\,.
\]
Putting these bounds together, we get
\[
ab\lesssim 1+(s+a)b\exp(-c_1bt^{-1}/2)\,.
\]
Since $\sup_{t>0}bt^{-1}\exp(-c_1bt^{-1}/2)\lesssim 1$, the following estimate holds
\[
ab\lesssim 1+st+ab\exp(-c_1bt^{-1}/2)\,.
\]
Recall that $a\le b$. Thus, an elementary bound $\sup_{t>0}b^2t^{-2}\exp(-c_1bt^{-1}/2)<\infty$ yields
\[
ab\exp(-c_1bt^{-1}/2)\le b^2\exp(-c_1bt^{-1}/2)\lesssim t^2\,.
\]
We finally get 
\[
ab\lesssim 1+st+t^2\lesssim 1+s^2+t^2
\]
and that proves the lemma.
\end{proof}
Now, given this lemma, we can argue in the following way. If $w,w^{-1}\in {\rm BMO}(\T)$,
then $w\in A_2(\T)$ and theorem~\ref{t2} yields
\begin{equation}\label{sade3}
\sup_{n}\int_\T |\phi_n|^{p}wd\theta<\infty\,, \quad  2\le p<p_{\rm cr}([w]_{A^2})\,.
\end{equation}
Therefore, for every $q\in [2,p)$, we can use H\"older's inequality
\begin{equation}\label{sade2}
\int_\T |\phi_n|^{q}d\theta=\int_\T |\phi_n|^{q}w^\beta w^{-\beta}d\theta\le 
\left(\int_\T |\phi_n|^{q\alpha}w^{\beta\alpha}d\theta\right)^{1/\alpha}\left( \int_\T w^{-\beta\alpha'}d\theta\right)^{1/\alpha'}
\end{equation}
and choose $\alpha\in (1,\infty)$ and $\beta>0$ such that $\beta\alpha=1, q\alpha=p$. The first factor in the right hand side of \eqref{sade2} is controlled by \eqref{sade3}. Since $w^{-1}\in {\rm BMO}(\T)$, the second factor is finite due to John-Nirenberg estimate and we get $\sup_{n}\|\phi_n\|_{L^{q}(\T)}<\infty$  as claimed in theorem \ref{dr}. This argument shows that theorem~\ref{t2} is qualitatively stronger than theorem \ref{dr}.\bigskip

\section{The Christoffel-Darboux Kernel and bounds for the associated Projection Operator}

In this section, we study the projection operators associated to $\{\phi_n (z,w)\}_{n \geq 0}$. Recall the Christoffel-Darboux kernel is defined as (see \cite{Simonbook}, p.120)
\[
	K_n (z,\zeta , w) = \sum\limits_{k=0}^n \phi_k (z,w) \overline{\phi_k (\zeta ,w )}.
\]
In particular, $K_n (z,\zeta, w)$ is integral kernel associated to the orthogonal projection operator $\proj{0}{n}{w}$ onto 
\mbox{$\mathrm{Span} \{\phi_0 , \ldots, \phi_n\}$} in $L^2_w(\T)$; see \cite{Simonbook} for more details. In this section, we prove that these projections are uniformly bounded:

\begin{Thm}\label{cd}
	Suppose $w \in A_2(\T)$, with $\gamma \ddd [w]_{A_2(\T)}$. Then, there exists $\epsilon_\gamma > 0$ such that 
	\[
		\sup\limits_{n} \|\proj{0}{n}{w}\|_{L^p_w(\T), L^p_w (\T)} < \infty
	\]
	for all $p \in [2 - \epsilon_\gamma ,2+ \epsilon_\gamma]$.
\end{Thm}

Recall (check \eqref{rev15}) that the Szeg\H{o} function $D$ can be introduced for any weight $w$ that satisfies $\log w\in L^1(\T)$.  We define the subspace $H_{2,w}(\T)$ as the closure of $\mathrm{Span} \{\phi_n\}_{n \geq 0}=\mathrm{Span} \{z^n\}_{n \geq 0}$ in $L_w^2(\T)$ metric. Denote  by  $\proj{0}{\infty}{w}$ the operator of  orthogonal projection  onto $H_{2,w}(\T)$ in $L^2_w(\T)$. By Beurling's theorem (\cite{Koos98}, p.79), function $f$ belongs to $H_{2,w}(\T)$ if and only if $f=D^{-1}g$ where $g$ is an element of the Hardy space $H_2(\T)$, e.g.,  $H_{2,w}(\T)=D^{-1}H_2(\T)$. Recall the standard notation that  $H_2(\T)$ denotes the restriction of functions in $H_2(\D)$ onto $\T$.
Since $w=|D|^2$, the map $g\rightarrow D^{-1}g$ is unitary isomorphism between $L^2(\T)$ and $L^2_w(\T)$.
The restriction of the same  map to $H^2(\T)$ is unitary isomorphism between $H_2(\T)$ and $H_{2,w}(\T)$.
 Finally, the orthogonal projection of $f\in L^2(\T)$ to  $H_2(\T)$ is given by $\lim_{r\to 1}\mathcal{C}(f,r\xi)$ (see \eqref{rev9} and \cite{vh}, p.2) where the limit exists for a.e. $\xi\in \T$. 
Thus, we can write
\begin{equation}\label{rev11}
	\proj{0}{\infty}{w} (f) (\xi) \ddd  \lim_{r\to 1}\frac{1}{D(\xi)} \mathcal{C}\Bigl(  f D,r\xi \Bigr), \quad \xi \in \T, 
\end{equation}
where $\mathcal{C}$ is Cauchy integral.

\begin{Lem}
	If $p\in (1,\infty)$ and $w^{1-p/2} \in A_{p} (\T)$, then $\proj{0}{\infty}{w}$ is bounded on $L^p_w(\T)$.
\end{Lem}
\begin{proof}
Let $\zeta\in \T$ and $z\in \D$. The Cauchy kernel in \eqref{rev9} can be written as
\[
\frac{1}{1-\bar \zeta z}=\frac 12\left(\frac{1+\bar \zeta z}{1-\bar\zeta z}+1\right)\,.
\]
The first term inside the parenthesis 
\[
\frac{1+\bar \zeta z}{1-\bar\zeta z}=\frac{\zeta+ z}{\zeta-z}
\]
is the so-called Schwarz kernel. 
Two real parts of Schwarz kernel is Poisson kernel \eqref{rev111} and its imaginary part, when restricted to $\T$, defines $\mathfrak{h}$ in \eqref{rev10}. Therefore, for $f\in L^p_w(\T)$, we can use \eqref{rev11} and \eqref{rev12} to get
\begin{eqnarray}\nonumber
|\proj{0}{\infty}{w} (f)|\lesssim \lim_{r\to 1}\frac{1}{|D|}\P(|fD|,r\xi)+\frac{1}{|D|}\int_\T |fD|d\theta+\left|\frac{1}{D}\mathfrak{h}\left(f D\right)\right|\\=
|f|+\frac{1}{|D|}\int_\T |fD|d\theta+\left|\frac{1}{D}\mathfrak{h}\left(f D\right)\right|\label{rev13}
\end{eqnarray}
due to (see p.11, \cite{Koos98}) and the identity
\[
\lim_{r\to 1}\P(g,r\xi)=g(\xi), \quad \text{a.e.}\,\, \xi\in \T
\]
which holds for $g\in L^1(\T)$. Since $f\in L^p_w(\T)$ and $w=|D|^2$, we get
\[
\left\|\frac{1}{|D|}\int_\T |fD|d\theta\right\|_{L^p_w(\T)}=
\left(\int_\T w^{1-p/2}d\theta\right)^{1/p}\cdot \left(\int_\T |f|\sqrt wd\theta\right)\,.
\]
Since $w^{1-p/2}\in A_p(\T)$ and $A_p(\T)\subset L^1(\T)$, the first integral converges. For the second one, we use H\"older's inequality
\[
\int_\T |f|\sqrt wd\theta=\int_\T (|f|w^{1/p})(w^{1/2-1/p})d\theta\le \left(\int_\T |f|^pwd\theta \right)^{1/p}\left(   \int_\T w^{(1/2-1/p)p'}d\theta\right)^{1/p'}\,.
\]
To show that the integral
\[
\int_\T w^{(1/2-1/p)p'}d\theta=\int_\T w^{\tfrac{(p-2)}{2(p-1)}}d\theta
\]
converges, we recall that $w^{1-p/2}\in A_p(\T)$ implies that 
$w^{\tfrac{(p-2)}{2(p-1)}}\in L^1(\T)$ as follows from the definition of $A_p(\T)$ given in \eqref{sd_00}.
We are left with estimating $L^p_w(\T)$ norm of the third term in \eqref{rev13}. The operator of harmonic conjugation $\mathfrak{h}$ is one of the basic singular integral operators and the Hunt-Muckenhoupt-Wheeden theorem claims (see, e.g., \cite{stein}, p.205) that  $\upsilon^{1/p}\mathfrak{h} \upsilon^{-1/p}$ is a bounded operator on $L^p(\T)$  if $\upsilon\in A_p(\T)$ and  $p\in (1,\infty)$. Since $w=|D|^2$ and $w^{1-p/2}\in A_p(\T)$, we get  statement of the lemma thanks to the formula
\[
\|w^{-1/2}\mathfrak{h}(w^{1/2}f)\|_{L_w^p(\T)}=\|w^{-1/2+1/p}\mathfrak{h}(w^{1/2-1/p}(w^{1/p}f))\|_{L^p(\T)}
\]
after one takes $\upsilon=w^{1-p/2}$ and notices that $\|w^{1/p}f\|_{L^p(\T)}=\|f\|_{L_w^p(\T)}$.
\end{proof}
This yields the following corollary.
\begin{Cor}\label{sad_bb}
	Let $w \in A_2(\T)$. Then, $\proj{0}{\infty}{w}$ is bounded on $L^p_w(\T)$ for all $p \in [4/3, 4]$. 
\end{Cor}
\begin{proof} The projection is self-adjoint operator in $L^2_w(\T)$. Therefore, by duality, it is enough to consider $p\in [2,4]$.
For $p=4$, we have $w^{-1}\in A_2(\T)\subset A_4(\T)$ and the previous lemma applies. If $p=2$, the projection operator has norm $1$. Thus, by Riesz-Thorin interpolation, we have an estimate for all $p\in [2,4]$.
\end{proof}

Define the projection operator onto $\mathrm{Span} \{\phi_n\}_{n \geq a+1}$ by 
\[
	\proj{a+1}{\infty}{w} \ddd \proj{0}{\infty}{w} - \proj{0}{a}{w}\, .
\]

When $w \in A_2(\T)$ and $p\in [4/3,4]$, $\{\proj{0}{n}{w}\}_{n \geq 0}$ is uniformly bounded on $L^p_w(\T)$ if and only if $\{\proj{n+1}{\infty,}{w} \}_{n \geq 0}$ is uniformly bounded on $L^p_w(\T)$. We will show the latter. To apply the same process as in section 3 for getting bounds for the polynomials $\{\phi_n\}$, one needs the following identities.
\begin{Lem}If $\proj{0}{n}{1}$ corresponds to the unperturbed case $w=1$, then
	\[
		\begin{cases}
	\proj{n+1}{\infty}{w} = (I - \proj{0}{n}{1}) \proj{0}{\infty}{w} + \proj{0}{n}{1} \proj{n+1}{\infty}{w} \\
			\proj{0}{n}{1} w \proj{n+1}{\infty}{w} = 0
		\end{cases} \, .
			\]
\end{Lem}
\begin{proof}
	To prove the first identity, first note that applying both operators to a function $f$ is the same as applying it to $\proj{0}{\infty}{w} f$, so it suffices to verify the identity for all functions in the range of $\proj{0}{\infty}{w}$ which is the closure of  finite sums $\sum_{j=0}^N a_j\phi_j(z)$. The formula then follows from $\proj{0}{n}{1} \phi_k = \phi_k$ for all $k \leq n$. To prove the second identity, it suffices to note that the range of $\proj{n+1}{\infty}{w}$ will be the closed span of $\{\phi_{n+1}, \phi_{n+2}, \ldots\}$; since $\phi_{n+j} \perp_{w} \{1,z, \ldots, z^n\}$, it follows that $\proj{0}{n}{1} w \phi_{n+j} = 0$ for all $j \geq 1$, whence the identity.
\end{proof}

{\it Proof of theorem \ref{cd}.}
By duality, it is sufficient to consider $p>2$.
Let $X_n \ddd w^{1/p} \proj{n+1}{\infty}{w} w^{-1/p}$ and 
$X_\infty\ddd w^{1/p}\cal{P}^w_{[0,\infty)}w^{-1/p}$. 
We need to estimate $\|X_n\|_{p, p}$. Rewriting the relations of the above lemma in terms of operators on $L^p (\T)$, we get
\[
	\begin{cases}
		X_n = w^{1/p}(I- \proj{0}{n}{1}) w^{-1/p} X_{\infty} + w^{1/p} \proj{0}{n}{1} w^{-1/p} X_n\\
		 w^{-1/p'} \proj{0}{n}{1} w^{1/p'} X_n = 0
	\end{cases}.
\]
Subtracting the bottom from the top and rearranging, we get back
\[
	(I - Q_{w, p}) X_n = w^{1/p} (I - \proj{0}{n}{1}) w^{-1/p} X_{\infty}.
\]
Notice that $\sup_n\|w^{1/p} (I - \proj{0}{n}{1}) w^{-1/p} X_{\infty}\|_{p,p}<\infty$ by Hunt-Muckenhoupt-Wheeden theorem and lemma~\ref{sad_bb}. Furthermore,  the proof of lemma \ref{l2} implies that $(I - Q_{w, p})$ on the left side of the equality has an inverse which is bounded in $L^p(\T)$ uniformly in $n$ for all $p \in [2, 2 + \epsilon_\gamma] \subseteq  [2,4]$ if $\epsilon_\gamma$ is small enough. Putting all of this together, we get
\[
	X_n = (I - Q_{w,p})^{-1} \Bigl(w^{1/p} (I - \proj{0}{n}{1}) w^{-1/p} X_{\infty}\Bigr) \, .
\]
Therefore, $\{X_n\}_{n \geq 0}$ is uniformly bounded, completing the proof.
\qed\smallskip

\section{Weights in $A_2(\T)$ and their Aleksandrov-Clark measures}

Several generalizations of $A_2(\T)$ and $A_\infty(\T)$ classes were studied in the literature (see, e.g., \cite{slavin}). We will need two definitions here.

\noindent{\bf Definition.}  We say that 
$w\in A_2^P(\T)$ if 
\begin{equation}\label{eq9}
[w]_{A_2^P(\T)}\ddd \sup_{z\in \D}\Bigl(\P(w,z)\P(w^{-1},z)\Bigr)<\infty
\end{equation}
and
$w\in A_\infty^P(\T)$  if
\begin{equation}\label{eq10}
[w]_{A_\infty^P(\T)}\ddd \sup_{z\in \D}\Bigl(\P(w,z)\exp(-\P(\log w,z)) \Bigr)<\infty\,.
\end{equation}
By Jensen's inequality, we have 
\begin{equation}\label{sd_93}
[w]_{A_\infty^P(\T)}\le [w]_{A_2^P(\T)}\,.
\end{equation}
The following lemma is part of the folklore of modern Harmonic Analysis, we  include its proof for completeness.
\begin{Lem}\label{ac} We have $A_2(\T)=A_2^P(\T)\subseteq A_\infty^P(\T)$.
\end{Lem}
\begin{proof}
By \eqref{sd_93}, we get the second inclusion. 
The inclusion $A_2^P(\T)\subseteq A_2(\T)$ follows from a bound
\[
\frac{1}{|I|^2}\left(\int_I wd\theta\right) \left(\int_I w^{-1}d\theta\right)\lesssim \P(w,z_I)\P(w^{-1},z_I)\,,
\]
where $z_I\ddd c_I(1-0.1|I|)$ and $c_I$ denotes the center of $I$. Thus, we 
only need to show $A_2(\T)\subseteq A_2^P(\T)$. Due to the rotational symmetry of $\D$, it is enough to take a point $z=1-\epsilon, \epsilon\in [0,1)$ and prove that
\begin{equation}\label{sd_yy}
\left(\int_{-\pi}^\pi \frac{\epsilon}{\epsilon^2+\theta^2}w(\theta)d\theta\right)\left(\int_{-\pi}^\pi \frac{\epsilon}{\epsilon^2+\theta^2}w^{-1}(\theta)d\theta\right)<C([w]_{A_2(\T)})\,.
\end{equation}
We can assume without loss of generality that 
\[
\langle w\rangle_{[0,\epsilon]}=1, \quad \langle w^{-1}\rangle_{[0,\epsilon]} \le [w]_{A_2(\T)}\,.
\]
In \cite{PeLe}, Lerner and Perez proved, in particular, that: 

{\it  Given $p\in (1,\infty)$, we have $w\in A_p(\R)$ if and only if
for every $\gamma>0$ there is $C(\gamma, [w]_{A_p})$ such that
\[
\frac{|E|}{|I|}\log^\gamma \left(\frac{|I|}{|E|}\right)\le C(\gamma, [w]_{A_p}) \left(\frac{w(E)}{w(I)}\right)^{1/p}\,,
\]
where $I$ is any interval in $\R$ and $E\subset I$.
}

Since each $w\in A_2(\T)$ can be considered as a $2\pi$-periodic weight on $\R$ with $[w]_{A_2(\R)}\lesssim [w]_{A_2(\T)}$, the result of Lerner and Perez holds for $\T$ as well. We take $p=2$, $E=[0,\epsilon], I=[0,x], 2 \epsilon<x<\pi$ to get
\[
\frac 1x\int_0^x w(s)ds\le C(\gamma, [w]_{A_2(\T)})\frac{x}{\epsilon} \log^{-2\gamma} \left(\frac{x}{\epsilon}\right)\,.
\]
Therefore when $\gamma > 1/2$ is fixed,
\begin{eqnarray*}
\int_0^\pi \frac{\epsilon w(x)}{\epsilon^2+x^2}dx\lesssim
	\epsilon^{-1} \int\limits_{0}^{2 \epsilon} w(x) dx +\epsilon\int_{2 \epsilon}^\pi \frac{w(x)}{x^2}dx\le
	C([w]_{A_2})+\epsilon \int_{2 \epsilon}^\pi \frac{1}{x^2}\left(\int_{2 \epsilon}^x w(\tau)d\tau  \right)'dx \lesssim
\\
	C([w]_{A_2})+\epsilon\int_{2 \epsilon} ^\pi w(x)dx+C(\gamma, [w]_{A_2(\T)})\int_{2 \epsilon}^\pi \frac{\log^{-2\gamma} (x/\epsilon)}{x}dx<C([w]_{A_2(\T)})\,,
\end{eqnarray*}
where in the second inequality we used that $A_2$ weights are doubling, along with our normalization. The integral over $[-\pi,0]$ can be estimated in the same way. Thus, 
\begin{equation}\label{sd_64}
\int_{\T}\frac{\epsilon w(x)}{\epsilon^2+x^2}dx<C([w]_{A_2(\T)})
\end{equation}
and we get a similar estimate for $w^{-1}$ because $w^{-1}\in A_2(\T)$.  We obtained \eqref{sd_yy} and the lemma is proved. 
\end{proof}

The following lemma was proved in \cite{BD2020} (see lemma 2 in this reference). We provide the sketch of the proof here.
\begin{Lem}
If $w\in A_\infty^P(\T)$ and $d\mu=\frac{w}{2\pi}d\theta$, then $\mu_{\alpha}$ is absolutely continuous and $d\mu_{\alpha}=\frac{w_{\alpha}}{2\pi}d\theta$ for every $\alpha\in \T$. Moreover, $w_{\alpha}\in A_\infty^P(\T)$.
\label{sad_ff}
\end{Lem}
\begin{proof}  Given probability measure $\mu:d\mu=\frac{w}{2\pi}d\theta+d\mus$, consider a generalized entropy 
$$
\K(\mu, z) = \log\P(\mu, z) - \P(\log  w, z), \qquad z \in \D.
$$
If we introduce $f$, the Schur function of measure $\mu$, through the formula
\begin{equation}\label{eq00}
\frac{1 + zf(z)}{1-zf(z)} =F(z)= \int_{\T}\frac{1 + \bar \xi z}{1 - \bar \xi z}\,d\mu(\xi), \quad z \in \D,\,\quad \xi=e^{i\theta}\,,
\end{equation}
then the straightforward but lengthy calculation  shows that
\begin{equation}\label{sad2}
\K(\mu, z) =\frac{1}{2\pi} \int_{\T}\log\left(\frac{1 - |zf(z)|^2}{1 - |f(\xi)|^2}\right)\,\frac{1-|z|^2}{|1 - \bar \xi z|^2}d\theta\,.
\end{equation}
On the other hand, it is known that the Schur function of each measure $\mu_\alpha$ is given by $f_\alpha=\alpha f$. Therefore, $\K(\mu_\alpha, z)=\K(\mu,z)$. Notice that $w\in A_\infty^P(\T)$ is equivalent to $\K(w,z)\in L^\infty(\D)$. Thus, if $w\in A_\infty^P(\T)$, then $\K(\mu_\alpha,z)\in L^\infty(\D)$. On the other hand, this condition implies that $\mu_\alpha$ has no singular part.
Indeed, if $d\mu_\alpha=\frac{w_\alpha}{2\pi} d\theta+d\mus^{(\alpha)}$ where  $\mus^{(\alpha)}$ is a singular measure, then  
\[
\log\left( \P(\mus^{(\alpha)},z)+\P(w_\alpha,z)\right)-\P(\log  w_\alpha,z)\le C, \quad z\in \D\,.
\]
This implies 
\[
\P(\mus^{(\alpha)},z)\le \P(\mus^{(\alpha)},z)+\P(w_\alpha,z)\le C\exp\left(\P(\log  w_\alpha,z)\right)\le C\P(w_\alpha,z)
\]
by Jensen inequality, hence, $\mus^{(\alpha)}=0$. 

\end{proof}

{\it Proof of theorem \ref{t4}.}
The first claim is immediate from lemma \ref{ac} and lemma \ref{sad_ff}.
Now, let us show that $w_\alpha\in A_2(\T)$. We will consider $w_{-1}=w_{\rm dual}$ only, the cases of other $\alpha$ can be handled similarly. We can write $F(e^{i\theta})=w+i\widetilde w$, where $\widetilde w$ is a harmonic conjugate function. Then, since $\Re F_{-1}=\Re F^{-1}=\Re F/|F|^2$, we get
\[
w_{\rm dual}=\frac{w}{w^2+\widetilde w^2}\,.
\]
Without loss of generality, we can consider an interval $I_\epsilon\ddd [-\epsilon, \epsilon]$ when checking $A_2(\T)$ condition for $w_{\rm dual}$.
We need to control
\begin{equation}\label{s1}
K\ddd \epsilon^{-2}\left(\int_{-\epsilon}^\epsilon \frac{w}{w^2+\widetilde w^2}d\theta\right)\left(\int_{-\epsilon}^\epsilon \frac{\widetilde w^2+w^2}{w}d\theta\right)
\end{equation}
under assumptions 
\begin{equation}\label{sd_k}
\langle w\rangle_{I_\epsilon}=1, \quad
\langle w^{-1}\rangle_{I_\epsilon}\le [w]_{A_2(\T)}\,.
\end{equation}
Clearly, 
\begin{equation}\label{s2}
\epsilon^{-2}\left(\int_{-\epsilon}^\epsilon \frac{w}{w^2+\widetilde w^2}d\theta\right)\left(\int_{-\epsilon}^\epsilon {w}d\theta\right)\lesssim  [w]_{A_2(\T)}
\end{equation}
by definition and we are left with estimating 
\begin{equation}\label{sd_s5}
\epsilon^{-2}\left(\int_{-\epsilon}^\epsilon \frac{w}{w^2+\widetilde w^2}d\theta\right)\left(\int_{-\epsilon}^\epsilon \frac{\widetilde w^2}{w}d\theta\right)\,.
\end{equation}
We can write 
\[
\widetilde w=h_1+h_2, \quad h_1\ddd \frak{h}({w\chi_{[-2\epsilon,2\epsilon]}}), \quad h_2\ddd \frak{h}(w\chi_{[-2\epsilon,2\epsilon]^c})\,,
\]
where $\frak{h}$ is harmonic conjugation, a standard singular integral operator.
Hence,
\[
\int_{-\epsilon}^\epsilon w^{-1}|h_1|^2d\theta\le
\int_{\T} w^{-1}|h_1|^2d\theta=\int_{\T} w^{-1}|\frak{h}(w^{1/2}\cdot w^{1/2}\chi_{[-2\epsilon,2\epsilon]}|^2d\theta\le C([w]_{A_2(\T)})\int_{-2\epsilon}^{2\epsilon}wd\theta
\]
if we use the Hunt-Muckenhoupt-Wheeden theorem with weight $w^{-1}\in A^2(\T)$ and $w^{-1/2}\frak{h}w^{1/2}$ applied to function $w^{1/2}\chi_{[-2\epsilon,2\epsilon]}$.
 In \eqref{sd_s5}, this gives the contribution
\begin{equation}\label{sd_s7}
\epsilon^{-2}\left(\int_{-\epsilon}^\epsilon \frac{w}{w^2+\widetilde w^2}d\theta\right)\left(\int_{-\epsilon}^\epsilon \frac{h_1^2}{w}d\theta\right)\le C([w]_{A_2(\T)})\epsilon^{-2}\left(\int_{-2\epsilon}^{2\epsilon}wd\theta\right) \left(\int_{-2\epsilon}^{2\epsilon}w^{-1}d\theta\right)\le C([w]_{A_2(\T)})\,.
\end{equation}
We are left with controlling
\begin{equation}\label{sd_34}
\epsilon^{-2}\left(\int_{-\epsilon}^\epsilon \frac{w}{w^2+\widetilde w^2}d\theta\right)\left(\int_{-\epsilon}^\epsilon \frac{h_2^2}{w}d\theta\right)\,.
\end{equation}
Notice that 
\[
h_2(\phi)=\Im U(e^{i\phi}), \quad \quad |\phi|<\epsilon\,,
\]
where
\[
U(\zeta)\ddd \frac{1}{2\pi}\int_{|\theta|>2\epsilon} \frac{e^{i\theta}+\zeta}{e^{i\theta}-\zeta}wd\theta, \quad \zeta\in \D\,.
\]
When $|\zeta-1|<\epsilon$, we have
\[
|U'(\zeta)|\lesssim \int_{|\theta|>2\epsilon} \frac{1}{|e^{i\theta}-1|^2}wd\theta\lesssim \epsilon^{-1}\int_\T \frac{\epsilon}{\theta^2+\epsilon^2}wd\theta\le \epsilon^{-1} C([w]_{A_2(\T)}),
\]
where we used the bound \eqref{sd_64}.
Therefore, 
\[
|\Im U(e^{i\phi})-\Im U(1-\epsilon)|\le C([w]_{A_2(\T)}), \quad |\phi|<\epsilon
\]
as follows from the Fundamental Theorem of Calculus. Therefore,
\begin{equation}\label{sd_ii}
\int_{-\epsilon}^\epsilon \frac{h_2^2}{w}d\theta\lesssim 
(\Im U(1-\epsilon))^2\int_{-\epsilon}^\epsilon w^{-1}d\theta+ C([w]_{A_2(\T)})\int_{-\epsilon}^\epsilon w^{-1}d\theta\,.
\end{equation}
The second term gives the following contribution in \eqref{sd_34}:
\begin{equation}\label{sd_341}
\epsilon^{-2}\left(\int_{-\epsilon}^\epsilon \frac{w}{w^2+\widetilde w^2}d\theta\right)C([w]_{A_2(\T)})\int_{-\epsilon}^\epsilon w^{-1}d\theta\le C([w]_{A_2(\T)}) \left(\langle w^{-1}\rangle_{I_\epsilon}\right)^2\le C([w]_{A_2(\T)})\,,
\end{equation}
where we used \eqref{sd_k}.
For the first term in \eqref{sd_ii}, recall that $\Re (F^{-1})= w/(w^2+\widetilde w^2)$ a.e. on $\T$ and estimate
\[
\epsilon^{-2}\left(\int_{-\epsilon}^\epsilon \frac{w}{w^2+\widetilde w^2}d\theta\right)(\Im U(1-\epsilon))^2\int_{-\epsilon}^\epsilon w^{-1}d\theta\lesssim \Bigl(\P(\Re (F^{-1}), 1-\epsilon)\cdot (\Im U(1-\epsilon))^2\Bigr)\cdot\Bigl( \epsilon^{-1}\int_{-\epsilon}^\epsilon w^{-1}d\theta\Bigr)\,.
\]
For the last factor, one can write
\[
\epsilon^{-1}\int_{-\epsilon}^\epsilon w^{-1}d\theta\lesssim [w]_{A_2(\T)}\,.
\]
Since $\Re (F^{-1})$ is harmonic, $\mu_{\rm dual}$ is absolutely continuous, and $\Re (F^{-1})=\Re F/|F|^2$, we get
\[
\P(\Re (F^{-1}), 1-\epsilon)\cdot(\Im U(1-\epsilon))^2=\frac{\Re F(1-\epsilon)}{|F(1-\epsilon)|^2}(\Im U(1-\epsilon))^2.
\]
Notice that our normalization gives
\begin{equation}\label{sd87}
1=(2\epsilon)^{-1}\int_{-\epsilon}^\epsilon wd\theta\lesssim \Re F(1-\epsilon)\sim 
 \int_{-\pi}^\pi \frac{\epsilon }{\theta^2+\epsilon^2}wd\theta\le C([w]_{A_2(\T)})\,,
\end{equation}
where the last bound is \eqref{sd_64}.
Let us compare $\Im U(1-\epsilon)$ and $\Im F(1-\epsilon)$. By definition of $F$ and $U$, 
\[
|U(1-\epsilon)-F(1-\epsilon)|\lesssim \frac{1}{\epsilon}\int_{-2\epsilon}^{2\epsilon} wd\theta\le C([w]_{A_2(\T)})\,.
\]
Thus, 
\begin{eqnarray*}
\frac{\Re F(1-\epsilon)}{|F(1-\epsilon)|^2}(\Im U(1-\epsilon))^2\lesssim  \frac{\Re F(1-\epsilon)}{|F(1-\epsilon)|^2}(|F(1-\epsilon)|^2+C([w]_{A_2(\T)}))\\<C([w]_{A_2(\T)}) \left(\Re F(1-\epsilon)+\frac{1}{\Re F(1-\epsilon)}\right)\,,
\end{eqnarray*}
which, thanks to \eqref{sd87}, is bounded by $C([w]_{A_2(\T)})$.
Summing up, we estimate $K$ in \eqref{s1} by $K\le C([w]_{A_2(\T)})$ and the lemma is proved.\qed

\section{Appendix: Fisher-Hartwig weights}

The Fisher-Hartwig weights are a large class of weights on the circle, which generalizes the class of Jacobi weights. It was at the focus of recent research (see, e.g., \cite{deift}) mainly due to some connections with probability and mathematical physics. 
For these weights, the asymptotics of polynomials is now well-understood \cite{deift}. In this section, we provide an upper bound for the function $p_{\rm cr}(t)$ using  some results obtained in \cite{Martinez06}. In particular, the analysis developed for Fisher-Hartwig weights will give us the proof of the following lemma.
\begin{Lem}If $t\in (1,2)$, we have $p_{\rm cr}(t)<C(t-1)^{-1/2}$.\label{sad_hh}
\end{Lem}
We provide its proof in the end of this section. For $\beta \geq 0$, consider the weight $w_{\beta} = |z-1|^{2 \beta}$ on the unit circle for  and the associated orthogonal polynomials $\{\Phi_n (z,w_{\beta})\}$. This is a particular choice for the Fisher-Hartwig weight with the single point of singularity located at $z=1$. Note that in order for $w_\beta \in A_2(\T)$, one needs $2 \beta < 1$, i.e. $\beta \in [0, \frac{1}{2})$. We start with the the following proposition:

\begin{Prop}
	Suppose $\beta \in [0, \frac{1}{2})$. Then 

	\[
[w_{\beta}]_{A_2(\T)} \sim \frac{1}{1- 4 \beta^2} \sim \frac{1}{1 - 2 \beta} \, .
		\]
	Furthermore, if $\beta\in [0,1/4]$, then
	\[
		[w_{\beta}]_{A_2(\T)} - 1 \sim \beta^2 \, .
\]
\end{Prop}
\noindent {\bf Remark.} The first asymptotics is useful in particular when $[w_{\beta}]_{A_2(\T)} > 2$, i.e.\ when our weight varies quite a bit, whereas when $[w_{\beta}]_{A_2(\T)} -1 < 1$, the second formula is more helpful. 

\begin{proof}  It is the straightforward calculation in which the integrals over intervals $I$ involved in the definition of $A_2(\T)$ can be explicitly computed and  estimated. We omit considering all cases here.
The formula which best explains the resulting bound is
\[
\langle \widetilde{w} \rangle_{I} \langle \widetilde{w}^{-1} \rangle_{I} = \frac{1}{1 - 4\beta^2}, \quad \widetilde w=|\theta|^{2\beta}
\] for $I = [0, a]$ and any $0\leq a\leq \pi$.

\end{proof}

The next proposition makes use of some statements from \cite{Martinez06}. Similar results for Jacobi weights were obtained in \cite{link1}.
\begin{Prop}
	Let $w_{\beta} = |z-1|^{2 \beta}$, $\beta \in [0,1/2)$. Then,
	\[
		\|\Phi_n (\cdot, w_{\beta})\|_{L^p_{w_\beta}(\T)} \sim_{\beta,p} \begin{cases} 1, \, &\text{ $2 \beta - p \beta +1 > 0 $} \\ \log  n, &\text{ $2 \beta - p \beta + 1 = 0$} \\ n^{-(2 \beta - p \beta +1)}, &\text{ $2 \beta - p \beta +1 < 0 $}\end{cases} .
			\]
			In particular, $\sup\limits_{n} \|\Phi_n (\cdot, w_{\beta})\|_{L^{p}_{w_{\beta}}(\T)} < \infty$ if and only if $p< 2+\frac{1}{\beta}$.
\end{Prop}

\begin{proof}
First, write
\[
	\|\Phi_n (\cdot, w_{\beta})\|_{L^p (w_{\beta})} ^p = \int_{|\theta| > \delta} |\Phi_n (z,w_{\beta})|^{p} w_{\beta}d\theta + \int_{|\theta| < \delta} |\Phi_n (z,w_{\beta})|^{p} w_{\beta}d\theta  \, ,
	\]
where $\delta$ is a parameter independent of $n$.
To control the first term, we use formula (1.13) of \cite{Martinez06} to get
	\[
\int_{|\theta| > \delta} |\Phi_n (z,w_{\beta})|^{p} w_{\beta}d\theta 		\le C(\beta,p,\delta) \int_{|\theta| > \delta} w_{\beta} ^{1- p/2}d\theta \le C(\beta,p,\delta)\,.
\]
As for the second term, using the asymptotics provided in (1.17) of \cite{Martinez06} and applying a change of variables $x = n \theta/2$, we get
	\[
		\int_{|\theta| < \delta} |\Phi_n (z,w_{\beta})|^{p} w_{\beta}d\theta \sim_\beta n^{p \beta - 2 \beta - 1} \int\limits_{0}^{\delta n/ 2} x^{2 \beta - p(\beta - 1/2)} |i J_{\beta + 1/2} (x) + J_{\beta - 1/2} (x)|^p dx \, ,
\]
where $J_{\nu} (x)$ is the Bessel function of the first kind.
One can then split this new integral in $x$ up into two: when $x\in (0,1)$ and when $x\ge 1$. We then use the known asymptotics for Bessel functions (see, e.g.,\ \cite{AS}) to get
\[
	\int_{|\theta| < \delta} |\Phi_n (z,w_{\beta})|^{p} w_{\beta}d\theta \sim_{\beta} n^{-(2\beta - p \beta + 1)} \Bigl(1 + \int\limits_{1}^{n \delta / 2} x^{2 \beta - p \beta} dx\Bigr) \sim_{\beta,p} \begin{cases} 1, \, &\text{ $2 \beta - p \beta +1 > 0 $} \\ \log n,  &\text{ $2 \beta - p \beta + 1 = 0$} \\ n^{-(2 \beta - p \beta +1)}, &\text{ $2 \beta - p \beta +1 < 0 $}\end{cases} \, .
	\]
In particular, this quantity is bounded precisely when $2 \beta - p \beta + 1  > 0$, i.e.\ when $\beta < \frac{1}{p-2}$. 
 The  proposition now follows from combining the given estimates.
\end{proof}
Now, we are ready to prove the main lemma of this section. 

{\it Proof of lemma \ref{sad_hh}.}
From the first proposition in appendix, we get $[w_\beta]_{A_2(\T)}-1\sim \beta^2$  if $\beta$ is small. The second proposition shows that $\sup_n \|\Phi_n(\xi,w_\beta)\|_{L^p_{w_\beta}(\T)}<\infty$ if and only if $p<2+\beta^{-1}$. Combining these results we get the statement of the lemma.
\qed

\bibliographystyle{plain} 
\bibliography{bibfile}

\end{document}